\documentclass{amsart}
\usepackage{rotating} 
\usepackage{circuitikz}
\usepackage{siunitx}
\usepackage{amsmath}
\usepackage{subfig}
\usepackage{array, makecell}
\usepackage{adjustbox}
\usepackage{multirow}
\usepackage[margin=1.3in]{geometry}
 \usepackage[foot]{amsaddr}

\newtheorem{problem}{Problem}

\begin{document}
\title[Optimal Control for Resistance-Type Boundary Conditions ]{An Optimal Control Approach to Determine Resistance-Type Boundary Conditions from in-vivo Data for Cardiovascular Simulations}
\author[Elisa Fevola et al.]{Elisa~Fevola$^1$, Francesco~Ballarin$^{2, 3}$, Laura~Jimenez-Juan$^4$, Stephen~Fremes$^5$, Stefano~Grivet-Talocia$^1$, Gianluigi~Rozza$^2$ and Piero~Triverio$^6$}
\address{$^1$ Department of Electronics and Telecommunications, Politecnico di Torino, Torino, Italy}
\address{$^2$ mathLab, Mathematics Area, SISSA, Trieste, Italy}
\address{$^3$ Department of Mathematics and Physics, Catholic University of the Sacred Heart, Brescia, Italy}
\address{$^4$ St Michael's Hospital and Sunnybrook Research Institute, University of Toronto, Toronto, Canada}
\address{$^5$ Sunnybrook Health Sciences Center and Sunnybrook Research Insititute, University of Toronto, Toronto, Canada}
\address{$^6$ Department of Electrical \& Computer Engineering, Institute of Biomedical Engineering, University of Toronto, Toronto, Canada}
\begin{abstract}
The choice of appropriate boundary conditions is a fundamental step in computational fluid dynamics (CFD) simulations of the cardiovascular system. 
Boundary conditions, in fact, highly affect the computed pressure and flow rates, and consequently haemodynamic indicators such as wall shear stress, which are of clinical interest.
Devising automated procedures for the selection of boundary conditions is vital to achieve repeatable simulations. 
However, the most common techniques do not automatically assimilate patient-specific data, relying instead on expensive and time-consuming manual tuning procedures. 
In this work, we propose a technique for the automated estimation of outlet boundary conditions based on optimal control. 
The values of resistive boundary conditions are set as control variables and optimized to match available patient-specific data. 
Experimental results on four aortic arches demonstrate that the proposed framework can assimilate 4D-Flow MRI data more accurately than two other common techniques based on Murray's law and Ohm's law.
\end{abstract}

\maketitle
\section{Introduction and motivation} \label{section:intro}
In recent years, large improvements have been made in developing patient-specific computational models to predict blood flow in patients affected by various cardiovascular diseases~\cite{taylor2009patient, gijsen2008strain}. The availability of data obtained from non-invasive medical imaging techniques, such as computed tomography (CT) and 4D Flow magnetic resonance imaging (MRI), combined with computational fluid dynamics (CFD), has led to more realistic blood flow modeling, with the possibility of investigating the origin and mechanisms behind different cardiovascular diseases.

\begin{figure}[t]
	\centering
	\includegraphics[scale=0.3]{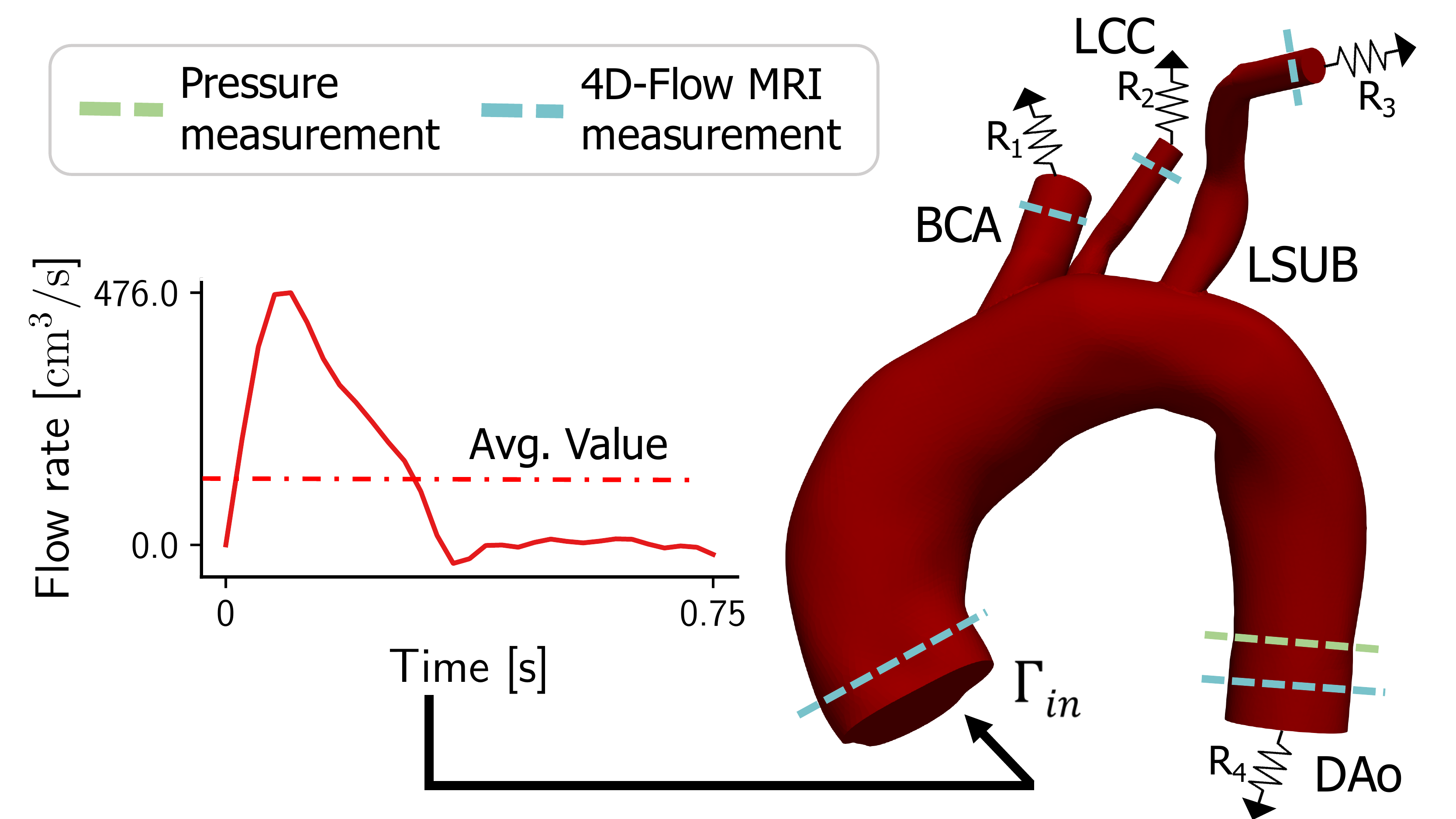}
	\caption{General configuration for the proposed optimal control approach on a patient's aortic arch. Patient-specific aortic arch haemodynamics is simulated by prescribing the measured inlet flow rate and the optimized resistive-like boundary conditions (BCs) at the outlets, being the outlets the brachiocephalic artery (BCA), the left common carotid artery (LCA), the left subclavian artery (LSUB) and the descending aorta (DAo). These BCs are calculated using pressure and flow measurements at several locations.}
	\label{fig:aorticarch}
\end{figure}

Mathematical models usually describe blood flow through Navier-Stokes equations, which for complex geometries, such as 3D models of the cardiovascular system, need to be solved numerically. Since it is unfeasible to discretize the entire cardiovascular system, Navier-Stokes equations are solved only on a portion of it, while the rest of the vasculature is accounted for by proper boundary conditions (BCs). Boundary conditions, in fact, model the upstream and downstream vasculatures not included in the 3D model by specifying the physiological conditions at the inlets and outlets of the computational domain of interest. An example of computational model for a human aortic arch can be found in Fig.~\ref{fig:aorticarch}, where the boundary conditions have to be imposed at the inlet $\Gamma_{in}$ and at the four outlets comprising the descending aorta (DAo) and the supra-aortic branches: brachiocephalic artery (BCA), left common carotid artery (LCC), and left subclavian artery (LSUB).
The selection of appropriate boundary conditions is always a crucial step when dealing with patient-specific cardiovascular models. 
Several studies in fact have shown how the choice of boundary conditions has a significant influence on predicted flow patterns as well as on clinically relevant parameters, such as wall shear stress (WSS)\cite{morbiducci2010outflow, pirola2017choice, van2011influence}.
This selection process often benefits from the availability of patient-specific in vivo measurements, so that BCs can be properly tuned to match available clinical data. These data might have been obtained with Phase Contrast Magnetic Resonance Imaging (PC-MRI), a technique for measuring blood velocity in specific sections of the vessel of interest, or with 4D-Flow MRI, a time-resolved PC-MRI with velocity encoding in all three spatial directions and 3D anatomic coverage~\cite{stankovic20144d}. 
For what concerns the inlet boundary condition, usually a velocity profile with a pulsatile waveform is prescribed. If patient-specific data are not available, physiological inlet flow waveforms are taken from the literature~\cite{pedley_1980}.
For what concerns outlet BCs, different strategies have been proposed, depending both on the complexity of the treated geometry and on the type of available clinical data.
A common choice consists in using constant pressure or traction boundary conditions, even if this greatly affects accuracy~\cite{pirola2017choice}.
A better solution resorts to the use of zero-dimensional (lumped parameter) models\cite{sankaranarayanan2005computational, sankaran2012patient, taylor2013computational, ge2018multi, grinberg2008outflow}, which, by prescribing a specific pressure-flow relationship at each outlet, lead to more physiological pressure and flow waveforms, and to a better approximation of in-vivo data. Lumped models may consist of single resistances or more sophisticated zero-dimensional models, such as the three-element Windkessel model~\cite{westerhof2009arterial}. 
The use of such outlet boundary conditions, however, leaves the user with the problem of finding appropriate numerical values for the lumped elements.
In this sense patient-specific in-vivo measurements, if available, can be used to guide the choice of the parameter values. 

The standard approach for cardiovascular CFD simulations consists in a first estimation of outlet resistive BCs through Murray's law~\cite{murray1926physiological}, which gives an indication of how blood flow splits in vessel bifurcations. These initial estimates are then refined to match available patient-specific data through a manual tuning process, until the user-desired accuracy is obtained \cite{boccadifuoco2018validation, mercuri2019tuning}.
Despite being straightforward to implement, this approach is sub-optimal, operator-dependent, time-consuming, and affected by non-repeatability.
To achieve robust and reliable patient-specific simulations, it is vital to rely on automated procedures for the selection of parameter values~\cite{bonfanti2019patient, spilker2010tuning}.
A first class of solutions is based on the use of a Kalman filter and its non-linear extensions, which performs the assimilation of measurements \textit{on the fly} and updates the unknown parameters accordingly~\cite{lal2017data, canuto2020ensemble}.
Arthurs et al.~\cite{arthurs2020flexible} proposed a reduced-order unscented Kalman filter for the sequential estimation of simple lumped parameter network parameters. Since the uncertainty estimation was limited to the unknown parameters, and not to the entire state space, the filtering operation was computationally tractable directly on 3D models. 
The data assimilation process, however, is still computationally expensive, as the estimation of \textit{n} parameters requires running \textit{n+1} forward simulations.

A different approach for parameter estimation relies on Bayesian inference, where the unknown parameters are described as random variables, with the goal of yielding confidence regions instead of point estimates and quantifying the uncertainty affecting simulation results~\cite{tran2017automated, schiavazzi2017patient, d2013uncertainty, perdikaris2016model}. Schiavazzi et al.~\cite{schiavazzi2017patient}, for example, used Bayesian estimation for automated tuning of 0D lumped model parameters to match clinical data, estimating also the uncertainty induced by errors in the measurements. Tran et al.~\cite{tran2017automated} proposed an iterative approach based on adaptive Markov chain Monte Carlo sampling to quantify confidence in local hemodynamic indicators, e.g., wall shear stress and oscillatory shear index (OSI).

Alternatively, variational approaches, such as optimal control, perform the parameter estimation by treating the governing equations as state system and the mismatch between the state solution and the patient-specific measurements as cost functional to be minimized \cite{koltukluouglu2018boundary, zainib2019reduced, funke2019variational, d2012applications}.
In practice, the unknown boundary conditions are chosen as control variables to minimize the cost functional.
As other advanced data assimilation methods, also optimal control suffers from a high computational cost, since it requires solving a constrained minimization problem. Nevertheless, recent studies have documented its successful application to several haemodynamics problems~\cite{quarteroni2003optimal, dede2007optimal, negri2015reduced}. For what concerns the calibration of hard-to-quantify boundary conditions, Tiago et al.\cite{tiago2017velocity}, for example, proposed a velocity control approach for the assimilation of velocity measurements in blood flow simulations. The framework was validated on an idealized 2D geometry and tested on a 3D model of a brain aneurysm, but with the assimilated velocity data still synthetically generated.
In the work by Koltukluouglu et al.~\cite{koltukluouglu2018boundary}, a data assimilation method based on optimal boundary control for 3D steady-state blood flow simulations was presented. The methodology was validated against real 4D-flow MRI data measured on a glass replica of a human aorta.
To reduce the computational cost of optimal control, Romarowski et al.~\cite{romarowski2018patient} identified a surrogate optimization problem for prescribing PC-MRI data as outlet boundary conditions, tuning the parameters of a three-element Windkessel model via a \textit{least-squares} approach. 
Zainib et al.~\cite{zainib2019reduced} proposed a reduced order framework for the application of optimal control to coronary artery bypass grafts, where synthetically-generated measurements were used to tune Neumann-type outlet boundary conditions, thus bringing optimal control closer to a real clinical setting.
Tiago et al.~\cite{tiago2014patient}, moreover, used an optimal control approach to address the uncertainty coming from the segmentation process, which affects the definition of the geometry and, in turn, wall shear stress and its derived measures.

In this work, we propose an automated method for the selection of resistive-type outlet boundary conditions based on the solution of an optimal control problem, which assimilates a set of patient-specific measurements, e.g., through 4D-Flow MRI imaging. The assimilation is performed by constraining the minimization of the objective functional to the solution of a steady Stokes problem. Since Stokes equations correspond to the linearized steady-state Navier-Stokes equations, they are less computationally expensive, while still providing an adequate estimation of boundary conditions. The control variables are the unknown resistances imposed at the outlet through the coupled multidomain method introduced by Vignon-Clementel et al.\cite{vignon2006outflow}. To ensure a better match between simulation results and \textit{in-vivo} data, the inlet flow waveform measured with 4D-Flow MRI is imposed as a patient-specific inlet boundary condition, as reported in Fig.~\ref{fig:aorticarch}. 
One of the novel aspects of this work is the use of more realistic boundary conditions with respect to Neumann and Dirichlet BCs, up to now the standard in an optimal control setting.
Its suitability for real clinical scenarios is demonstrated in this work by validating the framework on four aortic arches, reconstructed from medical images of real clinical cases. The presented method is used in these cases to set the outlet BCs on the descending aorta and the supra-aortic branches, assimilating real 4D-Flow MRI data. The boundary conditions obtained with the proposed method are compared to those obtained with two alternative calibration techniques, namely, Murray's law and Ohm's law, demonstrating the ability of optimal control to assimilate known physiological data consistently better. 
Moreover, an analysis of time-averaged wall shear stress and oscillatory shear index values obtained using the three different calibration methods is used to assess their effect on clinical haemodynamic indicators. 

This work is organized as follows: Section~\ref{section:methodology} presents the methodological details of the proposed optimal control approach. Results on four patient-specific anatomies are reported in Section~\ref{section:results}, while Section~\ref{section:limitations} analyses the advantages and limitations of the proposed approach and Section~\ref{section:conclusion} provides our conclusions and future perspectives.

\section{Methodology} \label{section:methodology}
\begin{figure}
	\centering
	\includegraphics[scale=0.25]{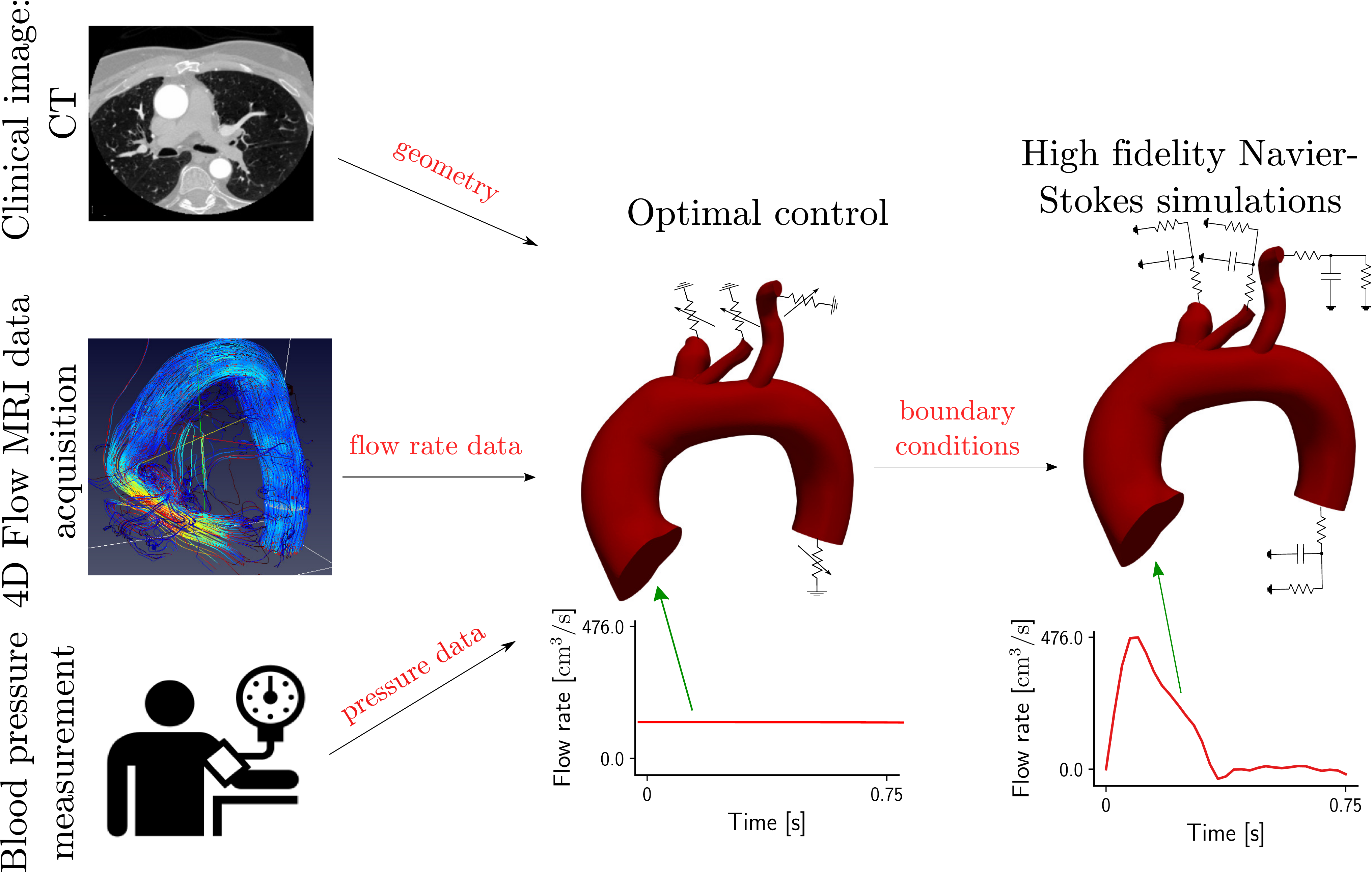}
	\caption{Scheme of the proposed framework for boundary condition estimation through optimal control.}
	\label{fig:pipeline}
\end{figure}

This section presents the methodological aspects of the proposed framework based on optimal control. 
Before solving the actual optimal control problem, some preliminary steps are required, namely, the extraction of patient-specific cardiovascular configurations from clinically acquired images, and the acquisition of patient-specific 4D-Flow MRI data. The proposed framework is represented in Fig.~\ref{fig:pipeline}. 

\subsection{Data acquisition and anatomical reconstruction} \label{geom_reconstruction}

The proposed methodology has been tested in four clinical cases from a single-center prospective study conducted at the Sunnybrook Health Sciences Centre in Toronto, Canada. The study was approved by the local ethics board and informed consent was obtained. Moreover, all the measurements were acquired non-invasively. Patients presented at the hospital for coronary bypass graft surgery. Between three and six weeks after surgery, a cardiac computer tomography (CT) was performed and anatomical information about their aorta and its supra-aortic branches was acquired using 320-detector row CT scanner (Aquilion One, Canon Medical Systems). From CT images, the vessels surface was reconstructed using the open-source package SimVascular~\cite{updegrove2017simvascular}.
The reconstructed volume was discretised into tetrahedral elements using TetGen~\cite{si2015tetgen}. After the CT scan, the blood velocity in the aorta and its branches was acquired in-vivo using a 4D-flow magnetic resonance imaging (MRI) sequence, using a  3T MRI scanner (MAGNETOM Prisma, Siemens Healthineers). The acquisition was performed using a 4D flow imaging sequence with retro-gating and adaptive navigator respiratory gating. Imaging parameters were as follows: velocity encoding=150 cm/s, field of view=200-420 mm x 248-368 mm, spatial resolution=1.9-3.5 x 2.0-3.2 x 1.8-3.5 mm$^3$, temporal resolution=39.9-47.2 ms, flip angle=8$^{\circ}$. After 4D-Flow MRI, diastolic blood pressure $P_{diast}$ and systolic blood pressure $P_{sys}$ were measured with brachial cuff-based method. The mean arterial pressure $P_{mean}$ was computed as
\begin{equation}\label{eq:pressure}
P_{mean} = \frac{P_{sys} + 2P_{diast}}{3}.
\end{equation}

\subsection{Determination of boundary conditions through optimal control}

In general, an optimal control problem consists in determining a number of control variables, which are unknown quantities, so that a cost functional is minimized, under the constraint of a system of equations describing the physical state of blood flow. 
The solution of this problem provides both an estimate for the unknown control variables, and the solution of the state variables (pressure and velocity). 
In the proposed framework, blood flow is modeled with the Stokes equations. This choice ensures a simpler optimal control problem with respect to its nonlinear version based on Navier-Stokes equations. 
Given the simplification that the Stokes equations introduce on the model, we will assess the effect on the quality of the determined output boundary conditions, by validating the obtained results on a Navier-Stokes model. This analysis is reported in Section~\ref{section:results}.

Considering the large diameter of the aorta and supra-aortic branches, we can safely assume that the blood behaves as a Newtonian fluid and that the viscosity is constant, since the dimension of blood particles is smaller compared to vessel diameters. Therefore, we adopt the incompressible Stokes equations as state equations modeling blood flow in the computational domain $\Omega$, which is the volume of the 3D aortic arch model reconstructed from medical images, as shown in Fig.~\ref{fig:aorticarch}. We refer to the boundaries of $\Omega$ as $\partial \Omega = \Gamma_{in} \cup \Gamma_w \cup \Gamma_{i}$, where $\Gamma_{in}$, $\Gamma_w$, and $\Gamma_{i}$ denote the inlet of the aorta, the vessel walls, and the outlets of the aortic arch, with $1\leq i \leq i_{max}$ ($i_{max}=4$ in the case of Fig.~\ref{fig:aorticarch}), respectively. The control variables are denoted as $R_i$.
Thus, state equations can be written in strong form as 
\begin{align} \label{eq:stokes}
\begin{cases}
-\nu \Delta \mathbf{v} + \nabla p = 0 &\text{in $\Omega$}, \\
\nabla \cdot \mathbf{v} = 0 &\text{in $\Omega$}, \\
\mathbf{v} = \mathbf{v}_{in} &\text{on $\Gamma_{in}$}, \\
\mathbf{v} = 0 &\text{on $\Gamma_w$}, \\
p = R_i\int_{\Gamma_i}\mathbf{v}\cdot \mathbf{n} \, d\Gamma_i &\text{on $\Gamma_i$, $1\leq i\leq i_{max}$,}
\end{cases}
\end{align}
where $\mathbf{v}$ is blood velocity, $p$ the pressure, $\nu=0.04$ dynes/cm$^2$s the dynamic viscosity, and $\mathbf{n}$ the outward normal to the outlets.
A plug profile is imposed at the inlet $\Gamma_{in}$ whose average value is extracted from the 4D-Flow MRI data. A no-slip condition is imposed at the vessel walls $\Gamma_w$, assumed to be rigid and non-permeable, while at the outlets $\Gamma_{i}$ a resistive type boundary condition is imposed, following the coupled multidomain method proposed by Vignon-Clementel et al~\cite{vignon2006outflow}.
Since the optimal control framework for partial differential equations~\cite{gunzburger2002perspectives} is typically derived starting from the weak formulation of the state equations, we introduce Hilbert spaces $V(\Omega)$ and $P(\Omega)$ for the velocity $\mathbf{v}$ and the pressure $p$, respectively. In particular, we choose $V(\Omega) = H^1(\Omega; \mathbb{R}^3)$ and $P(\Omega) = L^2(\Omega)$. Furthermore, we denote by $U = \mathbb{R}^{i_{max}}$ the space associated to the controls $\mathbf{R} = [R_1, \dots, R_{i_{max}}]^T \in U$. We note that $V(\Omega)$ and $P(\Omega)$ are function spaces (i.e., $\mathbf{v}(\mathbf{x})$ and $p(\mathbf{x})$ are functions depending on the spatial coordinate $\mathbf{x} \in \Omega$), while $U$ is simply an Euclidean space (i.e., $R_i$ are scalar numbers, and not functions). 

Starting from the strong form in~\eqref{eq:stokes}, the weak formulation is derived as: given $\mathbf{R} \in U$, find
$$\mathbf{v} \in V_{in}(\Omega) = \left\{\tilde{\mathbf{v}} \in V(\Omega): \tilde{\mathbf{v}}|_{\Gamma_{in}} = \mathbf{v}_{in} \text{ and } \tilde{\mathbf{v}}|_{\Gamma_{w}} = 0 \right\} \quad \text{and} \quad p \in P(\Omega)$$
such that
\begin{align} \label{eq:stokes_weak}
\begin{cases}
\displaystyle\nu \int_{\Omega}\nabla \mathbf{v} {\cdot} \nabla \mathbf{w}\, d\Omega - \int_{\Omega}p\,  (\nabla{\cdot}\mathbf{w})\, d\Omega + \sum_{i=1}^{i_{max}} R_i \displaystyle\int_{\Gamma_i}\mathbf{v}{\cdot} \mathbf{n}d\Gamma_i \displaystyle\int_{\Gamma_i} \mathbf{w}{\cdot}\mathbf{n} d\Gamma_i \\ \quad\quad + \sum_{i=1}^{i_{max}} \int_{\Gamma_i} \mathbf{w} {\cdot} \mathbf{n} (\mathbf{n} {\cdot} \nu \nabla \mathbf{v} \,\mathbf{n})d\Gamma_i - \sum_{i=1}^{i_{max}} \int_{\Gamma_i} \mathbf{w}{\cdot}\nabla \mathbf{v} \, \mathbf{n}  d\Gamma_i= 0 \quad &\text{in $\Omega$}, \\
\displaystyle\int_{\Omega} q\,(\nabla {\cdot}\mathbf{v})\, d\Omega = 0  &\text{in $\Omega$}, \\
\end{cases}
\end{align}
for every
$$\mathbf{w} \in V_0(\Omega) = \left\{\tilde{\mathbf{w}} \in V(\Omega): \tilde{\mathbf{w}}|_{\Gamma_{in}} = 0 \text{ and } \tilde{\mathbf{v}}|_{\Gamma_{w}} = 0 \right\} \quad \text{and} \quad q \in P(\Omega),$$
where $\mathbf{w}$ and $q$ are the test functions associated to velocity and pressure, respectively.

We now reformulate the term 
$$R_i \displaystyle\int_{\Gamma_i}\mathbf{v}{\cdot} \mathbf{n}d\Gamma_i \displaystyle\int_{\Gamma_i} \mathbf{w}{\cdot}\mathbf{n} d\Gamma_i$$
in $\eqref{eq:stokes_weak}_1$ to obtain an equivalent weak formulation which is more suitable for the forthcoming finite element discretization. We introduce a set of Lagrange multipliers $\lambda_i$, $1 \leq i \leq i_{max}$, defined as
$$\lambda_i = R_i \displaystyle\int_{\Gamma_i}\mathbf{v}{\cdot} \mathbf{n}d\Gamma_i.$$
We further denote by $\boldsymbol{\lambda} = [\lambda_1, \dots, \lambda_{i_{max}}]^T \in Z = \mathbb{R}^{i_{max}}$ the vector collecting the Lagrange multipliers, and $Z$ its associated space. Therefore, the equivalent weak formulation is: given $\mathbf{R} \in U$, find $\mathbf{v} \in V_{in}(\Omega),  p \in P(\Omega), \boldsymbol{\lambda} \in Z$ such that
\begin{align} \label{eq:stokes_weak_2}
\begin{cases}
\displaystyle\nu \int_{\Omega}\nabla \mathbf{v} {\cdot} \nabla \mathbf{w} d\Omega - \int_{\Omega}p\,(\nabla{\cdot}\mathbf{w}) d\Omega + \sum_{i=1}^{i_{max}} \int_{\Gamma_i} \lambda_i\, \mathbf{w} {\cdot} \mathbf{n} d\Gamma_i\\ \quad\quad + \sum_{i=1}^{i_{max}} \int_{\Gamma_i} \mathbf{w} {\cdot} \mathbf{n} \,(\mathbf{n}{\cdot}\nu \nabla \mathbf{v} \, \mathbf{n})d\Gamma_i -\sum_{i=1}^{i_{max}} \int_{\Gamma_i} \mathbf{w}{\cdot}\nabla \mathbf{v} \, \mathbf{n}  d\Gamma_i = 0 &\text{in $\Omega$}, \\
\displaystyle\int_{\Omega} q \, (\nabla{\cdot}\mathbf{v}) d\Omega = 0 &\text{in $\Omega$}, \\
\displaystyle\frac{1}{|\Gamma_i|}\int_{\Gamma_i} \lambda_i\, \eta_i\, d\Gamma_i - \int_{\Gamma_i} R_i\, \mathbf{v} \cdot \mathbf{n}\, \eta_i\, d\Gamma_i = 0\, \quad{\text{on $\Gamma_i$, $1\leq i \leq i_{max}$,}}  
\end{cases}
\end{align}
for every
$\mathbf{w} \in V_0(\Omega), q \in P(\Omega), \boldsymbol{\eta} = [\eta_1, \dots, \eta_{i_{max}}]^T \in Z$, where $\boldsymbol{\eta}$ collects the test ``functions'' associated to the Lagrange multipliers $\boldsymbol{\lambda}$.
This is the final system representing how the control $\mathbf{R}$ affects the underlying physics of the model. In order to set up the optimal control problem, a proper cost functional must also be defined.

We define the cost functional $J$ with the following form
\begin{equation}\label{eq:cost_functional}
\displaystyle J(\mathbf{v}, p) = \frac{\alpha_p}{2}\cdot \frac{\int_{\Gamma_p}||p-p_d||^2 d\Gamma_p}{\int_{\Gamma_p}||p_d||^2 d\Gamma_p} + \sum_{i=1}^{i_{max}} \frac{\alpha_i}{2} \cdot \frac{\left[ \int_{\Gamma_i}\mathbf{v}\cdot \mathbf{n} d\Gamma_i - Q_i\right]^2}{Q_i^2}.
\end{equation}

The first term in \eqref{eq:cost_functional} represents the normalized difference between the state pressure $p$ and the patient's average pressure $p_d$ measured at a specific cross section $\Gamma_p$; as explained in \ref{geom_reconstruction}, in this work $p_d$ was assumed equal to the mean arterial pressure, computed from the measured systolic and diastolic pressure. The second term, instead, represents the normalized difference between the calculated flow rate (obtained integrating velocity on the outlet section) and the flow rate $Q_i$ extracted from 4D-Flow MRI data at each outlet $\Gamma_i$. The choice of minimizing the normalized difference between simulated and measured quantities ensures that each term gives the same contribution to the optimization process, even when assimilating measurements with different orders of magnitude. The weights $\alpha_p$ and $\alpha_i$, however, can be used to change the contribution of each term individually. For the aortic arches under analysis, $\alpha_p$ and all $\alpha_i$ were set to 1, meaning that all measurements contribute equally to the minimization process.
The optimal control problem then reads: 
\begin{problem}
	Find $\mathbf{R}$ such that functional \eqref{eq:cost_functional} is minimized, under the constraint that $\mathbf{v}, p, \mathbf{\lambda}$ satisfy Equation \eqref{eq:stokes_weak_2}.
\end{problem}
To solve the optimal control problem, we adopt the adjoint-based Lagrangian approach~\cite{quarteroni2009numerical, hinze2008optimization, gunzburger2002perspectives, ito2008lagrange}. This approach models the optimal control problem as an unconstrained minimization problem, whose solution corresponds to the minimum of a properly defined \textit{Lagrangian functional}. In practice, the optimal solution is the one where all the derivatives of the Lagrangian functional vanish. 
The Lagrangian formulation requires the introduction of three new unknowns, the so-called adjoint variables, $\mathbf{z}\in V_0(\Omega)$, $b \in P(\Omega)$ and $\mathbf{t} = [t_1, \dots, t_{i_{max}}]^T \in Z$.
The Lagrangian functional for this problem takes the form
\begin{equation}\label{eq:lagrangian}
\begin{split}
\mathcal{L}(\mathbf{v}, p, \boldsymbol{\lambda}, \mathbf{R}, \mathbf{z}, b, \mathbf{t})&= J(\mathbf{v}, p) + \nu \int_{\Omega}\nabla \mathbf{v}{\cdot} \nabla \mathbf{z} d\Omega - \int_{\Omega} p\, (\nabla {\cdot}\mathbf{z})d\Omega + \sum_{i=1}^{i_{max}} \int_{\Gamma_i} \lambda_i \, \mathbf{z}{\cdot}\mathbf{n}\, d\Gamma_i \\ &
+ \sum_{i=1}^{i_{max}}\int_{\Gamma_i} \mathbf{z}{\cdot}\mathbf{n}\,(\mathbf{n} {\cdot}\nu \nabla \mathbf{v}\, \mathbf{n})d\Gamma_i  - \sum_{i=1}^{i_{max}} \int_{\Gamma_i} \mathbf{z} {\cdot} \nabla \mathbf{v}\, \mathbf{n}d\Gamma_i + \int_{\Omega}b \, (\nabla {\cdot}\mathbf{v})d\Omega \\ & + \sum_{i=1}^{i_{max}} \frac{1}{|\Gamma_i|} \int_{\Gamma_i}\lambda_i t_i d\Gamma_i - \sum_{i=1}^{i_{max}} \int_{\Gamma_i}R_i \mathbf{v} {\cdot} \mathbf{n} \, t_i d\Gamma_i.
\end{split}
\end{equation}
Given~\eqref{eq:lagrangian}, the optimality system can be obtained by imposing that the derivatives of $\mathcal{L}$ with respect to $(\mathbf{v}, p, \boldsymbol{\lambda}, \mathbf{R}, \mathbf{z}, b, \mathbf{t})$ must vanish, in short requiring $\nabla \mathcal{L} = 0$.
Taking the derivative of~\eqref{eq:lagrangian} with respect to $\mathbf{v}$ (denoted by $\mathcal{L}_{\mathbf{v}}$) in the direction $\mathbf{w}$ we obtain
\begin{subequations}
	\begin{flalign}\label{eq:lagrange_derivative_v}
	\langle\mathcal{L}_{\mathbf{v}}, \mathbf{w} \rangle &= \displaystyle\nu\int_{\Omega}\nabla \mathbf{w}{\cdot}\nabla \mathbf{z} d\Omega + \sum_{i=1}^{i_{max}} \mathbf{z}{\cdot} \mathbf{n}\,(\mathbf{n} {\cdot}\nu \nabla\mathbf{w}\, \mathbf{n})d\Gamma_i - \sum_{i=1}^{i_{max}} \int_{\Gamma_i} \mathbf{z}{\cdot} \nabla \mathbf{w}\, \mathbf{n}d\Gamma_i 
	+ \int_{\Omega}b \, (\nabla {\cdot}\mathbf{w}) d\Omega \\ &
	- \sum_{i=1}^{i_{max}} \int_{\Gamma_i}R_i \mathbf{w} {\cdot}\mathbf{n}\, t_i d\Gamma_i +\sum_{i=1}^{i_{max}} \frac{\alpha_i}{Q_i^2}\left[\int_{\Gamma_i}\mathbf{v}{\cdot}\mathbf{n} d\Gamma_i \int_{\Gamma_i}\mathbf{w}{\cdot} \mathbf{n} d\Gamma_i - Q_i\int_{\Gamma_i}\mathbf{w}{\cdot} \mathbf{n} d\Gamma_i\right ]=0, \notag &&
	\end{flalign}
	while taking the derivative of~\eqref{eq:lagrangian} with respect to $p$ in the direction $q$ we get
	\begin{flalign}
	\langle\mathcal{L}_p, q \rangle =\alpha_p\displaystyle\int_{\Gamma_p}(p-p_d)\, q \, d\Omega - \int_{\Omega} q \,(\nabla \cdot\mathbf{z}) d\Omega = 0. &&
	\end{flalign}
	Similarly, the derivatives of~\eqref{eq:lagrangian} with respect to $\lambda_i, R_i, \mathbf{z}, b, t_i$ are, respectively, 
	\begin{flalign}
	\langle\mathcal{L}_{\lambda_i}, m_i \rangle =\displaystyle \int_{\Gamma_i} m _i \mathbf{z} \cdot \mathbf{n} d\Gamma_i + \frac{1}{|\Gamma_i|}\int_{\Gamma_i}m_it_i d \Gamma_i = 0, \quad 1\leq i \leq i_{max},&&
	\end{flalign}
	\begin{flalign}
	\langle\mathcal{L}_{R_i}, r_i \rangle =- \int_{\Gamma_i}r_i \mathbf{v} \cdot \mathbf{n}\, t_i d\Gamma_i = 0,  &&
	\end{flalign}
	\begin{flalign}
\langle\mathcal{L}_{\mathbf{z}}, \mathbf{s} \rangle &= \displaystyle\nu \int_{\Omega}\nabla \mathbf{v} {\cdot} \nabla \mathbf{s} d\Omega - \int_{\Omega}p\, (\nabla {\cdot} \mathbf{s}) d\Omega + \sum_{i=1}^{i_{max}} \int_{\Gamma_i} \lambda_i \mathbf{s} {\cdot} \mathbf{n} d\Gamma_i + \sum_{i=1}^{i_{max}} \int_{\Gamma_i} \mathbf{s} {\cdot} \mathbf{n} (\mathbf{n}\nu \nabla \mathbf{v} {\cdot} \mathbf{n})d\Gamma_i\\ & - \sum_{i=1}^{i_{max}} \int_{\Gamma_i} \mathbf{s}{\cdot}\nabla \mathbf{v} \,\mathbf{n}  d\Gamma_i =0,\notag &&
\end{flalign}
	\begin{flalign}
	\langle\mathcal{L}_{b}, d\rangle = \displaystyle\int_{\Omega} d \,(\nabla\cdot \mathbf{v}) d\Omega =0, &&
	\end{flalign}
	\begin{flalign}
	\langle\mathcal{L}_{t_i}, \eta_i\rangle = \displaystyle\frac{1}{|\Gamma_i|}\int_{\Gamma_i} \lambda_i \eta_i d\Gamma_i - \int_{\Gamma_i} R_i \mathbf{v}\cdot \mathbf{n} \, \eta_i d\Gamma_i = 0, \quad 1\leq i \leq i_{max}.&&
	\end{flalign}
\end{subequations}
The presence in Eq.~\eqref{eq:lagrange_derivative_v} of a term containing the product of two integrals requires the use of an additional Lagrange multiplier.
We thus introduce the new variables $k_i = \int_{\Gamma_i} \mathbf{v}\cdot \mathbf{n} d\Gamma_i$, which we substitute in Eq.~\eqref{eq:lagrange_derivative_v}, and we add the following equations to the system
\begin{equation}\label{eq:lagrange_last}
\frac{1}{|\Gamma_i|}\int_{\Gamma_i}k_i \, c_i \, d\Gamma_i - \int_{\Gamma_i}\mathbf{v} \cdot \mathbf{n} \, c_i \, d\Gamma_i = 0,
\quad \forall c_i \in \mathbb{R},
\quad 1 \leq i \leq i_{max}.  
\end{equation}
Equations~\eqref{eq:lagrange_derivative_v} through~\eqref{eq:lagrange_last} form the so-called coupled optimality system, which we solve through a \textit{one-shot} approach~\cite{gunzburger2002perspectives, stoll2013all}, where the system  is solved directly for all the unknown variables.

We adopt an \textit{optimise-then-discretise} approach, meaning that we introduce the numerical discretisation of the problem after the derivation of the optimality system. In particular, we rely on Galerkin finite element methods for the solution of the discrete version of the system. 
The domain $\Omega$ was discretised into a finite mesh of size $h \in \mathbb{R}$, and so we introduce finite-dimensional solution spaces $V_h(\Omega), P_h(\Omega), U_h, Z_h$. In particular, we use Taylor-Hood elements for the velocity-pressure pair, i.e. $\mathbb{P}_2$ finite elements to define $V_h(\Omega)$ and $\mathbb{P}_1$ finite elements for $P_h(\Omega)$. Since the spaces associated to control and Lagrange multipliers are already finite-dimensional, we set $U_h = U$ and $Z_h = Z$. The discretised system is solved using the open-source libraries \textit{FEniCS}~\cite{alnaes2015fenics, logg2012automated} and \textit{multiphenics}~\cite{multiphenicswebsite}, the latter being an open-source library developed at SISSA mathLab for easy prototyping of problems characterized by a block structure and boundary restricted variables. The numerical solution of the problem is obtained by means of MUMPS~\cite{amestoy2001fully}, a parallel sparse direct solver.

\section{Numerical Results} \label{section:results}
\subsection{Synthetic data}
\begin{figure}[ht]
	\captionsetup[subfigure]{captionskip=7.5pt, margin=0cm,font=normalsize,format=plain,
		labelfont={bf,up},textfont={up}}
	\centering
	\subfloat[Convergence analysis of mean TAWSS.]{
		{\includegraphics[width=.45\textwidth]{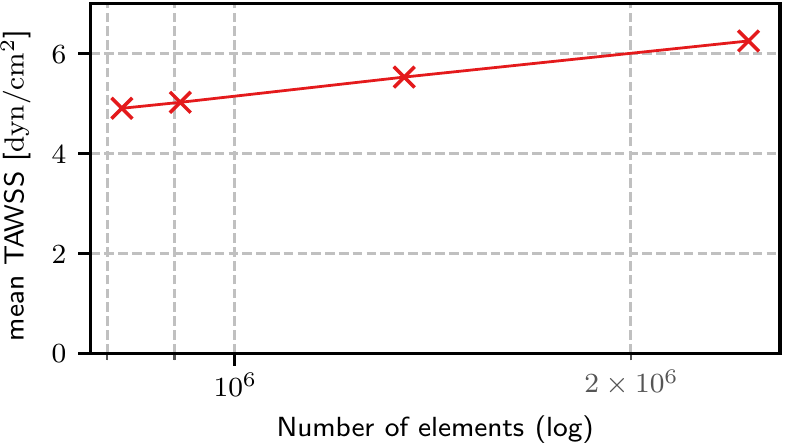}%
			\label{fig:wss-conv}}
	}
	\quad
	\subfloat[Convergence analysis of mean OSI.]{
		{\includegraphics[width=.45\textwidth]{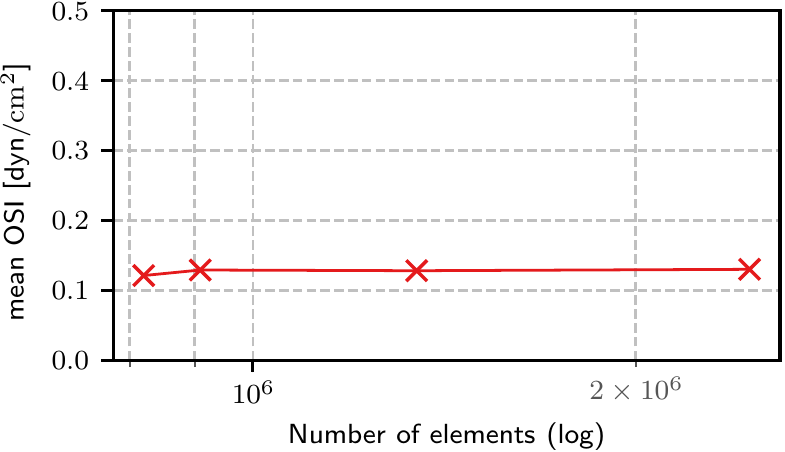}%
			\label{fig:osi-conv}}
	}
	\caption{Mesh convergence analysis.}
	\label{fig:mesh-convergence}
\end{figure}

A mesh convergence analysis was first conducted on the patient-specific meshes used for the experiments. In Figure~\ref{fig:mesh-convergence} the average TAWSS and OSI are plotted with respect to the number of elements in the mesh. The mesh used in the following experiments has about $2.5 \times 10^6$ elements, which has reached convergence for OSI, but not yet for WSS. This choice, however, is the best compromise in terms of accuracy and computational cost of the simulations.

\begin{figure}[ht]
	\captionsetup[subfigure]{captionskip=7.5pt, margin=0cm,font=normalsize,format=plain,
		labelfont={bf,up},textfont={up}}
	\centering
	\subfloat[Flow rate waveform at Brachiocephalic artery (BCA)]{
		{\includegraphics[width=.45\textwidth]{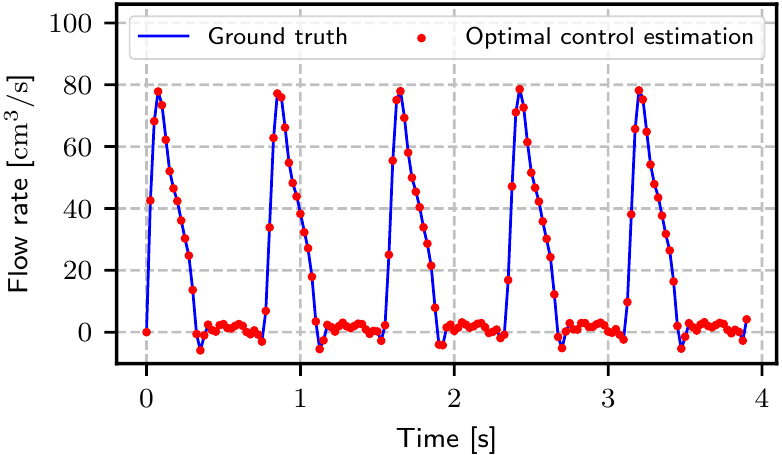}%
			\label{fig:synt-second}}
	}
	\quad
	\subfloat[Flow rate waveform at Left Common Carotid artery (LCC)]{
		{\includegraphics[width=.45\textwidth]{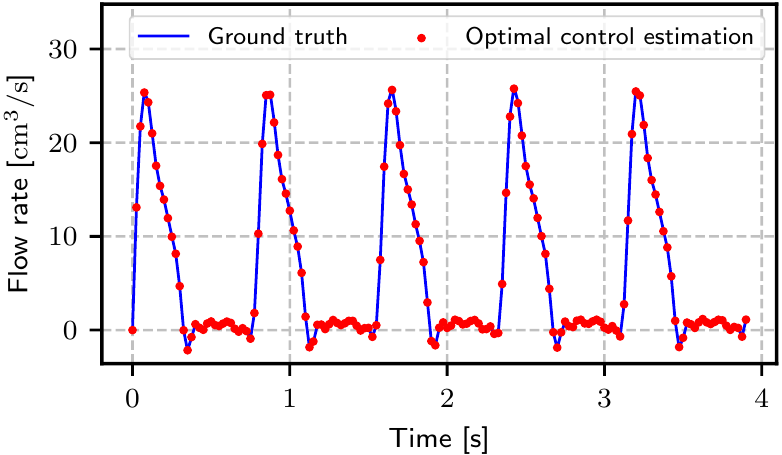}%
			\label{fig:synt-third}}
	}
	
	\subfloat[Flow rate waveform Left Subclavian artery (LSUB)]{
		{\includegraphics[width=.45\textwidth]{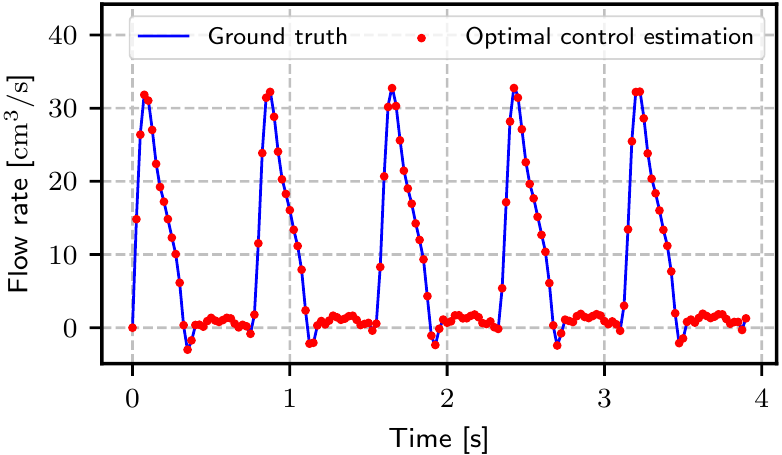}%
			\label{fig:synt-fourth}}
	}
	\quad
	\subfloat[Pressure waveform at Descending Aorta (DAo)]{
		{\includegraphics[width=.45\textwidth]{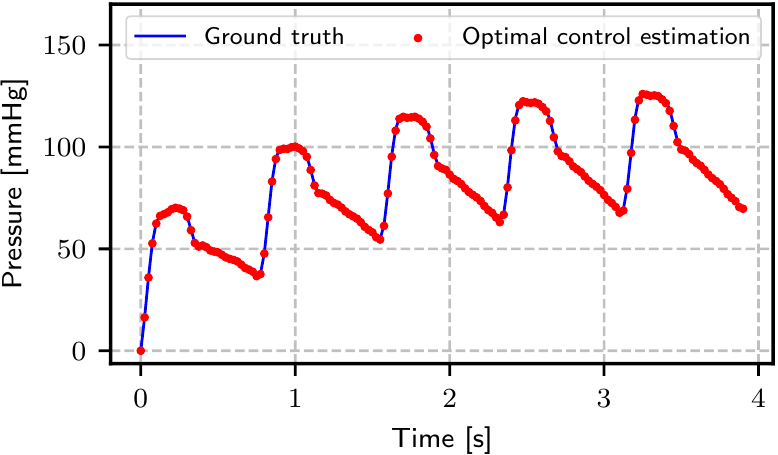}%
			\label{fig:synt-pressure}}
	}
	\caption{Comparison of flow rate waveforms at the three supra-aortic branches of case 1, and of pressure waveform at the Descending Aorta, for the test case with synthetically generated data.}
	\label{fig:synthetic_plots}
\end{figure}

The validity of the proposed approach was first tested on a test case (corresponding to the patient anatomy that will be labeled as ``case 1'' in the following) with synthetically generated data. 
To obtain the ground truth data, case 1 was simulated in SimVascular, imposing a velocity waveform at the inlet with average flow rate $Q =119.1 \text{ cm}^3\text{/s}$, using a blunt profile. At the outlets, three-elements Windkessel models were imposed. A physiological value of the total resistance at each outlet was arbitrarily chosen, and each resistance was then split into a proximal one, $R_{p,i} = 0.09R_i$, and a distal one, $R_{d,i} = 0.91R_i$, as suggested by Kim et al~\cite{kim2009coupling}. For the capacitance, a total value of $ 0.001$ cm$^5$/dyn was assumed~\cite{laskey1990estimation}, which was then split among the four outlets proportionally to their area. The distal pressure was set to 0 at all outlets. To reach periodic convergence, the simulation was run for five cardiac cycles, and the average flow rate at the four outlets and the average pressure were extracted from the last cardiac cycle. These data were fed into the optimal control tool presented in Section~\ref{section:methodology} to estimate the total outlet resistances. The estimated parameters were then used to set the outlet boundary conditions of a second simulation, again in Simvascular. Both resistances and capacitances were split adopting the same rules of the forward simulation.
Table~\ref{tab:synthetic} reports a comparison between original and estimated resistance values, together with original and estimated average flow rates at the outlets. Flow and pressure waveforms at the outlets are compared in Fig.~\ref{fig:synthetic_plots}, which shows how the resistances chosen with optimal control allow to reconstruct the original flow waveform with a relative error of 0.09\% for BCA, 0.09\% for LCC, 0.1\% for LSUB and 0.02\% for DAo. The pressure waveform is recovered with a relative error of 0.005\%.

\begin{table}[h!]
	\centering
	\begin{tabular}{|c|c|c|c|c|c|c|c|c|}
		\hline
		\multirow{2}{*}{} &  \multicolumn{4}{c|}{\textbf{Resistance (dyn$\cdot$s/cm$^5$)} } &  \multicolumn{4}{c|}{\textbf{Flow rate (cm$^3$/s) }} \\
		\cline{2-9}
		& \textbf{BCA} & \textbf{LCC} & \textbf{LSUB} & \textbf{DAo}  & \textbf{BCA} & \textbf{LCC} & \textbf{LSUB} & \textbf{DAo} \\
		\hline
		\textbf{Ground truth} & 7,000 & 21,000 & 16,000 & 1,700 & 19.02 & 6.30 & 8.13 & 80.02 \\
		\textbf{Estimated} & 6,937 & 20,846 & 15,991 & 1,685 & 19.05 & 6.30 & 8.09 & 80.03 \\
		\hline
	\end{tabular}
	\caption{Results of experiment with synthetic data.}
	\label{tab:synthetic}
\end{table}


\subsection{Patient-specific data}
Moving to the real cases with patient-specific data, the outlet boundary conditions estimated with the proposed optimal control framework were compared to those obtained with two alternative common techniques based on Murray's law and Ohm's law.
\paragraph{Murray's law}
Murray's law~\cite{murray1926physiological}, formulated by Cecil D. Murray in 1926, governs the branching pattern of vessels, such that the flow in each outlet is proportional to its cross-sectional area. In particular, the general form of Murray's law reads
\begin{equation}
\frac{Q_i}{Q_{tot}} = \frac{r_i^n}{\sum_i r_i^n}
\end{equation}
where $Q_i$ is the flow rate at the outlet, $r_i$ is the radius at the outlet, and $n$ for the aortic arch is conventionally set to 2.
Multiplying $r_i$ by $\pi$ to express the relationship in terms of outlet areas, and using $R \propto \frac{1}{Q}$, one can estimate the outlet resistances for the aorta and its main branches as
\begin{equation}\label{eq:murray}
R_i = \frac{\sum_j \left|\Gamma_j\right|}{\left|\Gamma_i\right|}R_{tot},
\end{equation}
where $\sum_j \left|\Gamma_j\right|$ is the sum of the area of all the aortic outlets, while $\left|\Gamma_i\right|$ is the area of the outlet to which resistance $R_i$ is associated. 
The total resistance $R_{tot}$ was computed as the mean pressure $p_d$ measured non-invasively on the patient, as reported in Section \ref{geom_reconstruction}, divided by the mean aortic flow rate $Q_0$ measured with 4D-Flow MRI, and then split among the outlet branches according to Equation~\eqref{eq:murray}. The application of Murray's law for predicting flow splitting has been largely investigated both on human and animal subjects by a number of studies~\cite{groen2010mri, trachet2011impact, cheng2007large,reneman2008wall}, which evidenced its validity on a large portion of the cardiovascular system, even if with some limitations on the first branches of the aortic arch~\cite{zamir1992relation}.
Note that, except for the average pressure and inlet flow rate, Murray's law does not take into consideration patient-specific measurements to compute flow splitting: instead, assuming that branches with larger cross-section have higher flow rates, it is solely based on the patient's anatomy. This feature justifies its use in those studies where in-vivo measurements of flow rates are not available~\cite{soulis2009influence, vincent2011blood, esmaily2011comparison}. 

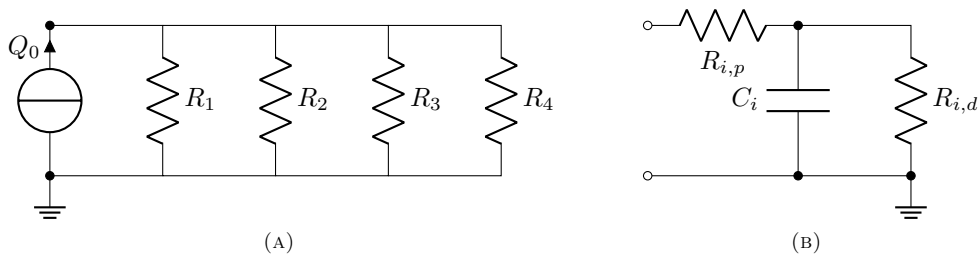
\begin{figure}[t]
	\centering
	\subfloat[]{
		\begin{circuitikz} \draw
			(0,0) node[ground] {}
			to[I, i=$Q_0$, *-*] (0,2) -- (1.5,2)
			to[R, l=$R_1$] (1.5,0) -- (0,0)
			(1.5,2) -- (3,2)
			to[R, l=$R_2$] (3,0) -- (1.5,0)
			(3,2) -- (4.5,2)
			to[R, l=$R_3$] (4.5,0) -- (3,0)
			(4.5,2) -- (6,2)
			to[R, l=$R_4$] (6,0) -- (4.5,0)
			;\end{circuitikz}  
		\label{fig:circuit}}%
	\qquad
	\subfloat[]{
		\centering
		\begin{circuitikz} \draw
			(0,0) node[anchor=east]{}
			to[short, o-*] (2,0)
			to[C, l=$C_i$] (2,2)
			to[R, l=$R_{i,p}$, *-o] (0,2)  
			(2,2) --(3.5, 2)
			to[R, l=$R_{i,d}$, -*] (3.5, 0) node[ground] {}
			(3.5, 0) -- (2, 0)
			;\end{circuitikz}   
		\label{fig:3wk}}
	\caption{Adopted equivalent circuits. On the left, equivalent circuit used for Ohm's law method. On the right, 3-element Windkessel model used as outlet boundary condition for unsteady Navier-Stokes simulations.}
\end{figure}

\paragraph{Ohm's law}
Based on the Ohm's law, it is possible to set outlet resistances by taking advantage of the analogy between the cardiovascular system and electrical circuits. In particular, knowing the mean pressure $p_d$ computed from diastolic and systolic pressure measured non-invasively and the outlet flow rates $Q_{i}$ measured with 4D-Flow MRI, outlet resistances can be computed as
\begin{equation}\label{eq:ohm}
R_i = \frac{p_d}{Q_i}.
\end{equation}
This method is also common~\cite{pirola2017choice} and it is based on the idea of performing the parameter estimation on 0D models~\cite{pant2014methodological} but, differently from Murray's law, it requires the availability of outlet flow rates measured in-vivo. The measured flow rates usually respect the mass conservation principle with approximately a $15\%$ deviation~\cite{dyverfeldt20154d}, due to measurement uncertainty. For this reason, we developed a minimization problem which, using Ohm's law, estimates outlet resistances $R_i$ while trying to impose the mass conservation principle. In this way, we try to compensate for the intrinsic inconsistency in the data, which is not accounted for in Ohm's law of Eq.~\eqref{eq:ohm}.
This is achieved by approximating the aorta with its equivalent 0D model, which is represented in Fig.~\ref{fig:circuit}, and minimizing the cost function
\begin{equation}\label{eq:ohm_fminsearch}
J_{ohm} = \alpha_p\frac{||R_{tot}\cdot Q_0-p_d||^2}{||p_d||^2} + \sum_{i=1}
^{i_{max}}\alpha_i \frac{||\frac{p_d}{R_i}-Q_i||^2}{||Q_i||^2}.
\end{equation}
The cost function was minimized using the Matlab function \textit{fminsearch}, which employs a derivative-free simplex search method~\cite{lagarias1998convergence}. The cost functional reported in Equation~\eqref{eq:ohm_fminsearch} deliberately replicates the one used inside the optimal control framework, reported in Equation~\eqref{eq:cost_functional}, with the first term representing the normalized difference between the computed pressure and the measured one ($p_d$), and the second the outlet flow rates $Q_i$. However, optimal control relies on Stokes equations as a 3D-model of the underlying system, whereas here a 0D approximation is used. 

\begin{figure}
	\centering
	\includegraphics[width=.9\linewidth]{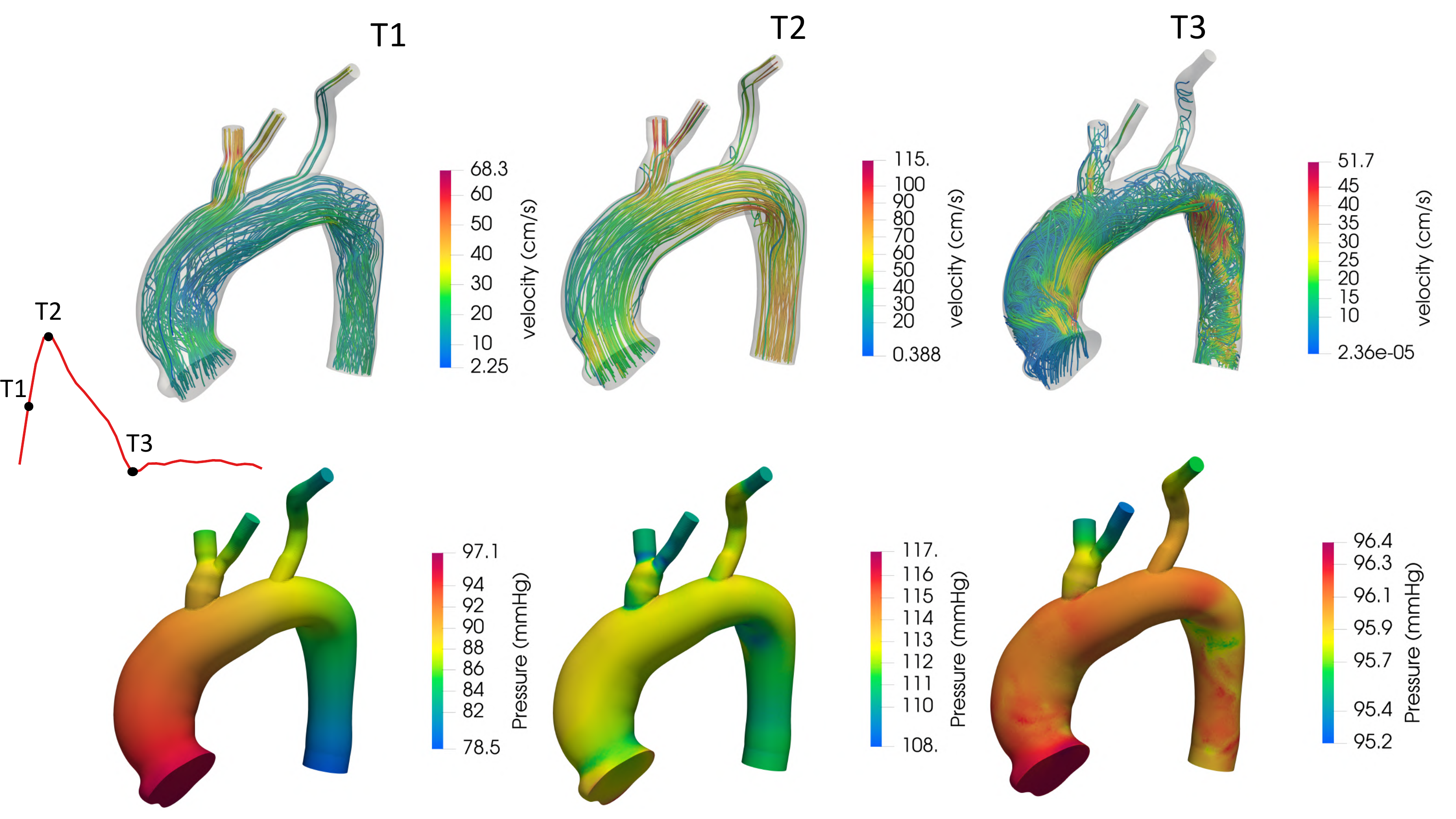}
	\caption{Velocity streamlines (top row) and pressure distributions (bottom row) for case 3 at three different time instants. The time points T1 (mid-systolic acceleration, left column), T2 (peak systole, middle column) and T3 (diastole, right column) are defined along a flow waveform shown on the left.}
	\label{fig:velocity_pressure}
\end{figure}
The experiments described in this section were conducted on four patient anatomies, obtained as described in Section~\ref{geom_reconstruction}. 
Table~\ref{tab:net_flows} reports for each case the measured flow rates. Specifically, the column labeled \textit{Net flow} indicates the difference between measured inlet flow and outlet flows, thus quantifies the violation of mass conservation on the measurements. 
The inlet flow reported in Table~\ref{tab:net_flows} was obtained by subtracting 4\% to the flow measured in the ascending aorta, which corresponds to the total coronary circulation~\cite{eslami2020effect}, not considered in the models.
The presence of flow rate inconsistencies in 4D-Flow MRI data could be due to the finite resolution of 4D-flow MRI, motion artifacts, and the presence of noise, especially in presence of complex helical and vortical flows~\cite{contijoch20204d}. 
The amount of net flow is below $15\%$ for all the cases under analysis, which is considered acceptable for 4D-Flow MRI measurements according to Dyverfeldt et al~\cite{dyverfeldt20154d}. While in case 1 and case 4 the net flow is significant, for cases 2 and 3 it is practically negligible.

\subsection{Boundary conditions estimation: Stokes model}\label{stokes_results}
The main results are reported in Table~\ref{table:ocp} where, for each case, the first row reports the patient's average pressure, measured non-invasively after MRI, and the flow rates measured in-vivo with 4D-Flow MRI at the four outlets: BCA, LCC, LSUB and DAo. 
The values obtained solving Stokes equations on the 3D geometry, with resistance values calculated by means of Ohm's law (Eq.~\eqref{eq:ohm}) are reported in the second row, while the third row contains the results obtained with the optimization based on Ohm's law as in Eq.~\eqref{eq:ohm_fminsearch}. The fourth contains the results for Murray's Law. Finally, the fifth row shows the results of the proposed methodology based on optimal control.
In all simulations, the inlet flow was imposed using a plug profile.
For case 2 and case 3, where the net flow is negligible, both Ohm's law and optimal control properly assimilate the available data, with a relative error of less than $1\%$ on all the outlet flow rates. For case 1 and case 4, where the net flow is high, optimal control outperforms the other techniques, achieving the smallest relative errors on pressure and BCA, LCC, LSUB flow rates. In presence of a large net flow, we cannot expect a perfect assimilation of measurements, as they are intrinsically non-physical. In that case, optimal control shows a smaller sensitivity to inconsistencies in the data, leading to a physical solution which is closer to measurements with respect to the other techniques.
As expected, the optimization based on Ohm's law is more accurate than Ohm's law for case 1 and 4, while the two methods are basically equivalent for case 2 and 3. 
As already pointed out earlier, as Murray's law is the only technique which does not take into account the measured flow rates to set outlet boundary conditions, the corresponding solution is the one that most deviates from patient measurements.
The solution of the optimal control problem required an average of 6.75 minutes (wall clock CPU time), running on 18 Lenovo SD530 nodes, each with 40 Intel "Skylake" cores and 202 GB RAM.

\begin{table}[h!]
	\centering
	\begin{tabular}{|c|c|c|c|c|}
		\hline
		\makecell{\textbf{Case} \\ \textbf{number} }& \makecell{\textbf{Inlet flow} \\ \textbf{(cm$^3$/s)}} & \makecell{\textbf{Outlet flow}\\ \textbf{(cm$^3$/s)}} &\makecell{ \textbf{Net flow} \\ \textbf{(cm$^3$/s)}} & \makecell{ \textbf{Mass conservation} \\ \textbf{violation (\%)}} \\
		\hline
		\textbf{1} & 119.10 & 103.46 & 15.64 & 13\% \\
		\textbf{2} & 107.00 & 107.18 & -0.18 & 0.17\% \\
		\textbf{3} & 125.63 & 125.60 & 0.03 & 0.02\% \\
		\textbf{4} & 103.00 & 90.21 & 12.79 & 12.4 \% \\
		\hline
	\end{tabular}
	\caption{Table summarizing, for each case, the measured inlet flow rate, the sum of the measured outlet flow rates, the difference between the two (Net flow), and the percentage of mass conservation violation.}
	\label{tab:net_flows}
\end{table}

	\begin{table}
	\centering
	\begin{adjustbox}{max width = \textwidth}
		\begin{tabular}{|c|c|c|c|c|c|c|}
			\hline
			\multirow{2}{*}{\makecell{\textbf{Case} \\ \textbf{number}}} & \multirow{2}{*}{\textbf{Method}} & \multirow{2}{*}{\makecell{\textbf{Pressure} \\ \textbf{(mmHg)}}} & \multicolumn{4}{c|}{\textbf{Flow rate (cm$^3$/s) (\% error w.r.t. measurements)}} \\
			\cline{4-7}
			&       &                     &  \textbf{BCA}  & \textbf{LCC}  & \textbf{LSUB} & \textbf{DAo}  \\
			\hline
			\multirow{5}{*}{\textbf{5}}&\textbf{Measurements}& 98.7 & 15.9 & 5.98 & 8.48 & 73.1\\
			& \textbf{Ohm's law}  & 113 (15\%) & 18.3 (15.3\%) & 6.89 (15\%) & 9.77 (15\%) & 84.2 (15\%) \\
			& \textbf{Opt. based on Ohm's law} & 108 (9.5\%)& 17.6 (11\%) & 6.55 (9.5\%) & 9.28 (9.4\%) & 85.8 (17\%) \\
			& \textbf{Murray's law} & 98.8 (0.1\%) & 19.2 (20\%) & 6.19 (3.5\%) & 7.48 (-11\%) & 86.2 (18\%) \\
			&\textbf{Proposed} & 98.7 (0\%) & 16.6 (4.4\%) & 6.08 (1.7\%) & 8.67 (2.2\%) & 87.8 (20\%) \\
			\hline
			\multirow{5}{*}{\textbf{9}}&\textbf{Measurements}& 105& 13.2 & 6.71 & 7.47 & 79.8 \\
			& \textbf{Ohm's law} & 105 (0\%) & 13.1 (-0.7\%) & 6.68 (-0.4\%) & 7.44 (-0.4\%) & 79.5 (-0.4\%)\\
			& \textbf{Opt. based on Ohm's law} & 105 (0\%) & 13.1 (-0.7\%) & 6.70 (-0.14\%) & 7.45 (-0.26\%) & 79.6 (-0.3\%) \\
			& \textbf{Murray's law} & 106 (0.9\%)  & 17.0 (29\%) & 6.28 (6\%) & 9.30 (24\%) & 74.2 (-7\%) \\
			&\textbf{Proposed} & 106 (0.9\%) & 13.2 (0\%) & 6.71 (0\%) & 7.47 (0\%) & 79.5 (-0.4\%) \\
			\hline
			\multirow{5}{*}{\textbf{10}}&\textbf{Measurements}& 103  & 19.0 & 11.3 & 10.5 & 84.8  \\
			& \textbf{Ohm's Law} & 103 (0\%) & 19.0 (0\%) & 11.3 (0\%) & 10.5 (0\%) & 84.9 (0.1\%) \\
			& \textbf{Opt. based on Ohm's law} & 104 (0.65\%) & 19.1 (0.5\%) & 11.4 (0.62\%) & 10.5 (0\%) & 84.9 (0.1\%)  \\
			& \textbf{Murray's law} & 104 (0.65\%) & 18.1 (-5\%) & 9.37 (-17\%) & 7.22 (-31\%) & 91.2 (7\%)  \\
			&\textbf{Proposed} & 103 (0\%) & 19.0 (0\%) & 11.3 (0\%) & 10.5 (0\%) & 85.0 (0.2\%)  \\
			\hline
			\multirow{5}{*}{\textbf{11}}&\textbf{Measurements}& 100  & 9.87 & 4.32 & 6.92 & 69.1 \\
			& \textbf{Ohm's law}& 115 (15\%)  & 11.3 (14\%) & 4.94 (14\%) & 7.91 (14\%) & 78.9 (14\%) \\
			& \textbf{Opt. based on Ohm's law}& 109 (9\%) & 10.8 (9\%) & 4.68 (8.3\%) & 7.49 (8.2\%) & 80.1 (16\%)\\
			& \textbf{Murray's law} & 101 (1\%) & 15.3 (55\%) & 4.09 (-5\%) & 5.09 (-26\%) & 78.6 (13.7\%) \\
			&\textbf{Proposed} & 100 (0\%)& 10.1 (2.6\%) & 4.37 (1.1\%) & 7.04 (1.7\%) & 81.5 (18\%) \\
			\hline
		\end{tabular}
	\end{adjustbox}
	\caption{Comparison of pressure and flow rates for Stokes simulations with boundary conditions obtained with the different methods.}\label{table:ocp}
	
\end{table}

\begin{table}[h!]
	\centering
	\begin{tabular}{|c|c|c|c|c|c|}
		\hline
		\multirow{2}{*}{\makecell{\textbf{Case} \\ \textbf{number}}} & \multirow{2}{*}{\textbf{Method}} &  \multicolumn{4}{c|}{\textbf{Resistance (dyn$\cdot$s/cm$^5$) }} \\
		\cline{3-6}
		&  & \textbf{BCA} & \textbf{LCC} & \textbf{LSUB} & \textbf{DAo} \\
		\hline
		\multirow{4}{*}{\textbf{1}} & \textbf{Ohm's law} & 8,288 & 21,979 & 15,518 & 1,800 \\
		& \textbf{Opt. based on Ohm's law} & 8,190 & 21,981 & 15,511 & 1,679 \\
		& \textbf{Murray's law} & 6,837 & 21,242 & 17,591 & 1,527 \\
		& \textbf{Proposed} & 7,941 & 21,609 & 15,153 & 1,497 \\
		\hline
		\multirow{4}{*}{\textbf{2}} & \textbf{Ohm's law} & 10,690 & 21,003 & 18,847 & 1,764 \\ 
		& \textbf{Opt. based on Ohm's law} & 10,737 & 20,987 & 18,867 & 1,767  \\
		& \textbf{Murray's law} & 8,306 & 22,451 & 15,161 & 1,900 \\
		& \textbf{Proposed} & 10,699 & 21,006 & 18,864 & 1,773 \\
		\hline
		\multirow{4}{*}{\textbf{3}} & \textbf{Ohm's law} & 7,248 & 12,142 & 13,094 & 1,624 \\
		& \textbf{Opt. based on Ohm's law} & 7,255 & 12,145 & 13,090 & 1,629 \\
		& \textbf{Murray's law} & 7,646 & 14,723 & 19,102 & 1,513 \\
		& \textbf{Proposed} & 7,249 & 12,131 & 13,078 & 1,621 \\
		\hline
		\multirow{4}{*}{\textbf{4}} & \textbf{Ohm's law} & 13,600 & 31,060 & 19,391 & 1,943 \\
		& \textbf{Opt. based on Ohm's law} & 13,500 & 31,069 & 19,399 & 1,815 \\
		& \textbf{Murray's law} & 8,751 & 32,800 & 26,358 & 1,708 \\
		& \textbf{Proposed} & 13,166 & 30,500 & 18,903 & 1,634 \\
		\hline
	\end{tabular}
	\caption{Resistance values chosen by the different methods.}
	\label{tab:resistances}
\end{table}

\subsection{Boundary conditions estimation for inaccurate measurements}
In presence of a large net flow, the presented approach compensates for the inaccuracy in the data by adjusting the outlet boundary conditions. This means that the inconsistencies in the measurements are resolved entirely at the outlets, while the measurement imposed at the inlet is considered deterministic. It would be desirable, instead, to gain some flexibility in the assimilation of the inlet flow, which is equally affected by uncertainty. In this case, a modification of the approach presented in Section~\ref{section:methodology} is possible, which estimates both the inlet and outlet BCs by means of an additional control at the inlet Dirichlet boundary condition. In particular, the velocity at the inlet in~\ref{eq:stokes} is expressed as
\begin{equation}
\mathbf{v} = u_{in} \cdot \mathbf{v_{in}} \quad\text{on $\Gamma_{in}$},
\end{equation}
where $u_{in}$ is a scalar control variable and $v_{in}$ is the inlet velocity profile, with average value equal to the one measured with 4D-Flow MRI.
By choosing the best value for $u_{in}$, the optimal control problem will be able to change the inlet flow rate to better assimilate the available data. This obviously requires a slight modification of the cost functional, with an additional term for assimilating the flow rate $Q_{in}$ at the inlet $\Gamma_{in}$:
\begin{equation}\label{eq:cost_functional_inlet}
\displaystyle J(\mathbf{v}, p) = \frac{\alpha_p}{2}\cdot \frac{\int_{\Gamma_p}||p-p_d||^2 d\Gamma_p}{\int_{\Gamma_p}||p_d||^2 d\Gamma_p} + \sum_{i=1}^{i_{max}} \frac{\alpha_i}{2} \cdot \frac{\left[ \int_{\Gamma_i}\mathbf{v}\cdot \mathbf{n} d\Gamma_i - Q_i\right]^2}{Q_i^2}  + \frac{\alpha_{in}}{2} \cdot \frac{\left[ \int_{\Gamma_{in}}\mathbf{v}\cdot \mathbf{n} d\Gamma_{in} - Q_{in}\right]^2}{Q_{in}^2}.
\end{equation}
Moreover, the Dirichlet control requires weakly imposing the inlet condition by means of Lagrange multiplier, thus increasing the final size of the system of equations~\eqref{eq:lagrange_derivative_v} - \eqref{eq:lagrange_last}.
In table~\ref{table:ocp_new} we report the results with this alternative formulation for cases 1 and 4, which had the largest Net flow. As expected, the net flow is now resolved acting on all boundary conditions, leading to lower differences between measured and simulated flow rates.
It is worth mentioning that the additional control leads to an increase in the dimensions of the problem, and consequently to larger computational costs (an average of 10 minutes of wall clock CPU time, running on 24 Lenovo SD530 nodes, each with 40 Intel "Skylake" cores and 202 GB RAM).

\begin{figure}[ht]
	\captionsetup[subfigure]{captionskip=7.5pt, margin=0cm,font=normalsize,format=plain,
		labelfont={bf,up},textfont={up}}
	\centering
	\subfloat[case 1]{
		{\includegraphics[width=.45\textwidth]{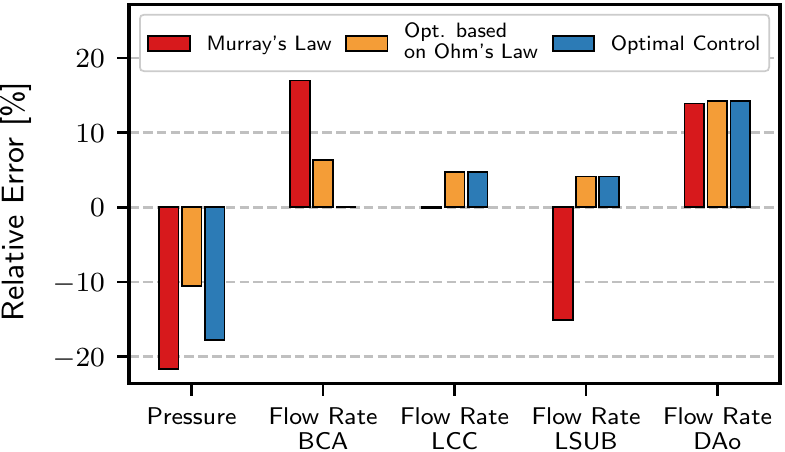}%
			\label{fig:sub-first}}
	}
	\quad
	\subfloat[case 2]{
		{\includegraphics[width=.45\textwidth]{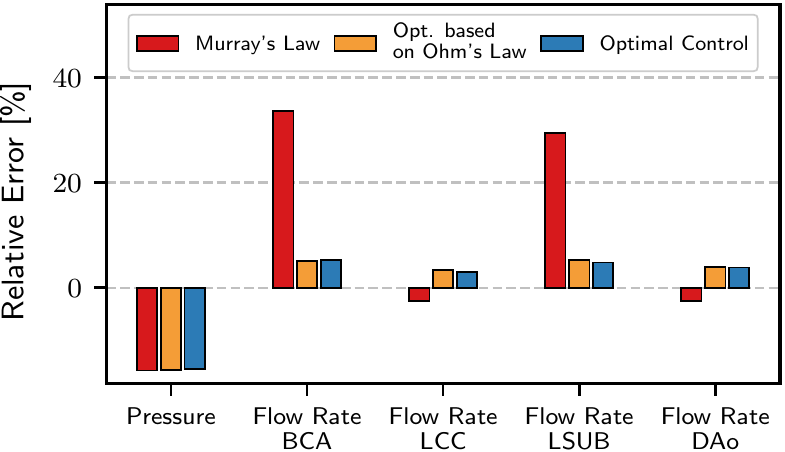}%
			\label{fig:sub-second}}
	}

	\subfloat[case 3]{
		{\includegraphics[width=.45\textwidth]{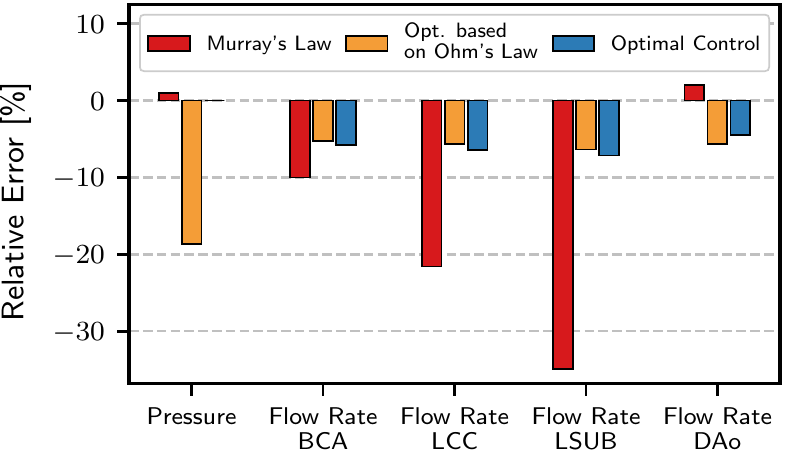}%
			\label{fig:sub-third}}
	}
	\quad
	\subfloat[case 4]{
		{\includegraphics[width=.45\textwidth]{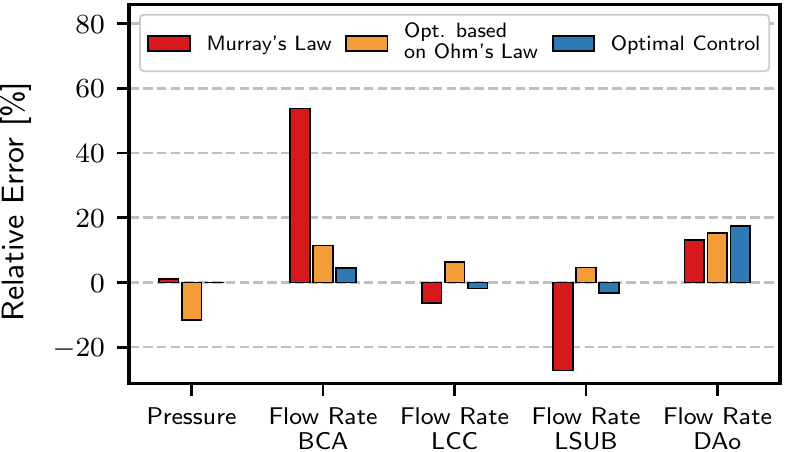}%
			\label{fig:sub-fourth}}
	}
	\caption{Comparison of pressure and outlet flow rates for time-dependent Navier-Stokes simulations. The histograms report the relative difference with respect to the corresponding patient measurements.}
	\label{fig:histograms}
\end{figure}

\begin{table}[h!]
	\begin{center}
	\begin{adjustbox}{max width = \textwidth}
		\begin{tabular}{|c|c|c|c|c|c|c|c|}
			\hline
			\multirow{2}{*}{\makecell{\textbf{Case} \\ \textbf{number}}} & \multirow{2}{*}{\textbf{Method}} & \multirow{2}{*}{\makecell{\textbf{Pressure} \\ \textbf{(mmHg)}}} & \multicolumn{5}{c|}{\textbf{Flow rate (cm$^3$/s) (\% error w.r.t. measurements)}} \\
			\cline{4-8}
			&       &                 &\textbf{Inlet}    &  \textbf{BCA}  & \textbf{LCC}  & \textbf{LSUB} & \textbf{DAo}  \\
			\hline
			\multirow{2}{*}{\textbf{1}}&\textbf{Measurements}& 98.7 & 119.1 & 15.9 & 5.98 & 8.48 & 73.1\\
			&\textbf{Proposed} & 98.7 (0\%) & 107.9 (-9\%) & 16.1 (1.2\%) & 6.01 (0.5\%) & 8.53 (0.6\%) & 77.3 (5.7\%) \\
			\hline
			\multirow{2}{*}{\textbf{4}}&\textbf{Measurements}& 100  & 103.0 & 9.87 & 4.32 & 6.92 & 69.1 \\
			&\textbf{Proposed} & 100 (0\%)& 94.3 (-8\%) & 9.95 (0.8\%) & 4.33 (0.23\%) & 6.96 (0.6\%) & 73.0 (6\%) \\
			\hline
		\end{tabular}
		\end{adjustbox}

		\caption{Comparison of measurements and pressure and flow rates for Stokes simulations obtained with the extended optimal control formulation, which controls both the inlet and the outlets boundary conditions.}
		\label{table:ocp_new}
	\end{center}
\end{table}

\subsection{Unsteady Navier-Stokes simulations}
An accurate representation of blood flow in the aorta is generally obtained through unsteady Navier-Stokes simulations, which provide a more realistic and accurate time evolution of blood flow. Assuming laminar regime, the incompressible Navier-Stokes equations take the form
\begin{align} \label{eq:navierstokes}
\begin{cases}
-\nu \Delta \mathbf{v} + \nabla p + \frac{\partial \mathbf{v}}{\partial t} + (\mathbf{v} \cdot \nabla) \mathbf{v} = 0 &\text{in $\Omega$}, \\
\nabla \cdot \mathbf{v} = 0 &\text{in $\Omega$.}
\end{cases}
\end{align}

We used the resistance values $R_i$ reported in Table~\ref{tab:resistances} as outlet boundary conditions of high-fidelity unsteady Navier-Stokes simulations performed using SimVascular~\cite{updegrove2017simvascular}. The purpose of this analysis is twofold. First, it verifies that the flow splitting obtained with a steady linear Stokes model is still valid in a non-linear, unsteady scenario. Second, it allows to analyse the impact that boundary conditions obtained with different estimation techniques have on wall shear stress-related indicators, which are clinically relevant. 
For unsteady Navier-Stokes simulation, the type of outlet boundary condition which best replicates realistic flow conditions is the three-element Windkessel model~\cite{pirola2017choice} represented in Fig.~\ref{fig:3wk}.
Referring to Fig.~\ref{fig:3wk}, each resistance $R_i$ was split into a proximal one, $R_{p,i} = 0.09R_i$, and a distal one, $R_{d,i} = 0.91R_i$, as suggested by Kim et al~\cite{kim2009coupling}. For the capacitance, a total value of $ 0.001$ cm$^5$/dyn was assumed~\cite{laskey1990estimation}, which was then split among the four outlets proportionally to their area. The distal pressure was set to 0 at all outlets. The dynamic viscosity was set to $\nu=0.04$ dynes/cm$^2$s, a rigid wall model was assumed and a plug profile was imposed at the inlet. The inlet waveform was the one extracted from 4D flow MRI. The time-step value for the transient simulations was set to 0.5 ms, and 5 cardiac cycles were simulated. The results reported here refer to the last cardiac cycle. Each simulation required an average of 8 hours (clock wall CPU time), running on 4 Lenovo SD530 nodes, each with 40 Intel "Skylake" cores and 202 GB RAM. 

Figure~\ref{fig:velocity_pressure} reports the velocity streamlines and pressure distribution for case 3 at three different time instants along the cardiac cycle.
A comparison of pressure and outlet flow rates obtained with Navier-Stokes simulations using outlet boundary conditions estimated with the three different techniques introduced previously is reported in Fig.~\ref{fig:histograms}. The histograms represent, for each outlet, the relative difference of the average flow rate with respect to the measured one using Murray's law, Ohm's law, or the proposed method. Also for this simulations, the inlet flow was imposed using a plug profile. As expected, the errors of the obtained flows with respect to the measured ones increased when moving from a Stokes model to a Navier-Stokes one, mostly due to the non-linearity and time-dependency introduced by the latter. Nevertheless, similar trends to those of the Stokes experiments reported in Table~\ref{table:ocp} can be observed when moving to Navier-Stokes simulations. In particular, Murray's law remains the method providing the largest deviations from measured flow rates, while both Ohm's law and optimal control are closer to assimilated data. With the exception of case 3, optimal control is still the method which best replicates measured flow rates. These results show that the BCs estimated with a linear, steady Stokes model proved to be a good choice when moving to high-fidelity Navier-Stokes simulations, thus supporting the approach of estimating BCs on a linearised Stokes model.

\begin{figure}[ht]
	\captionsetup[subfigure]{captionskip=7.5pt, margin=0cm,font=normalsize,format=plain,
		labelfont={bf,up},textfont={up}}
	\centering
	\subfloat[Flow rate waveform at Brachiocephalic artery (BCA)]{
		{\includegraphics[width=.45\textwidth]{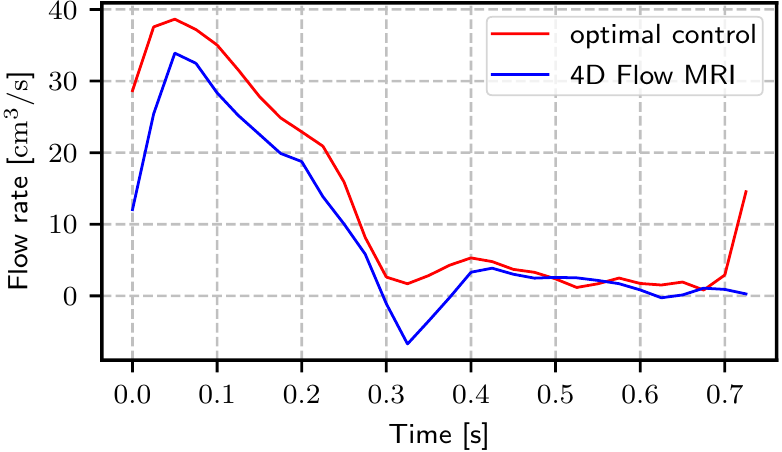}%
			\label{fig:4flow-bca}}
	}
	\quad
	\subfloat[Flow rate waveform at Left Common Carotid artery (LCC)]{
		{\includegraphics[width=.45\textwidth]{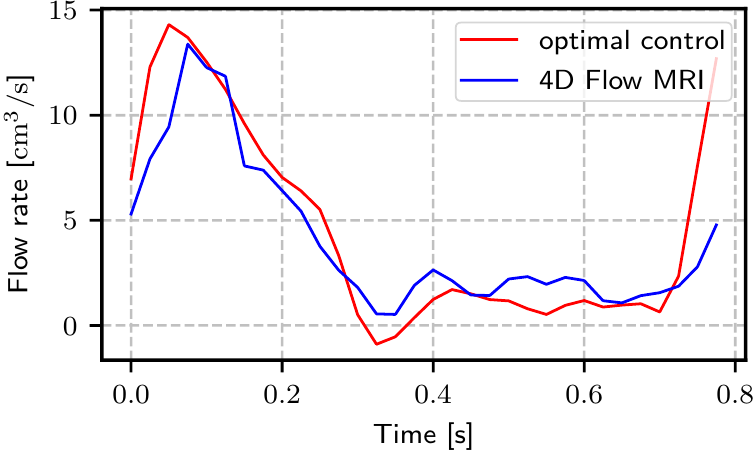}%
			\label{fig:4flow-lcc}}
	}

	\subfloat[Flow rate waveform Left Subclavian artery (LSUB)]{
		{\includegraphics[width=.45\textwidth]{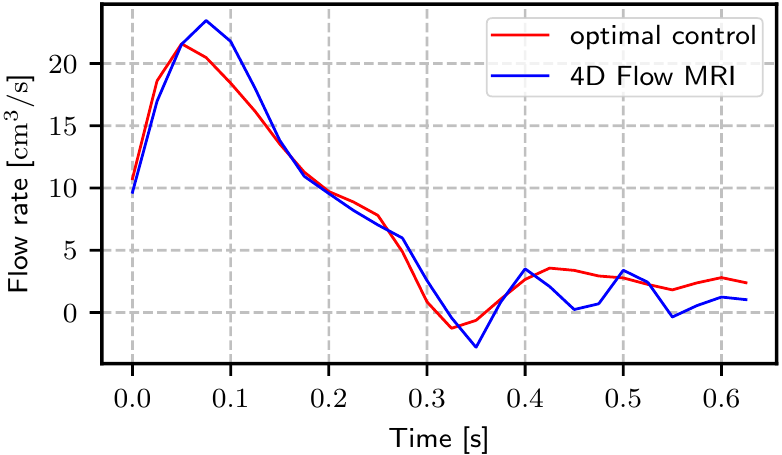}%
			\label{fig:4flow-lsub}}
	}
	\quad
	\subfloat[Flow rate waveform at Descending Aorta (DAo)]{
		{\includegraphics[width=.45\textwidth]{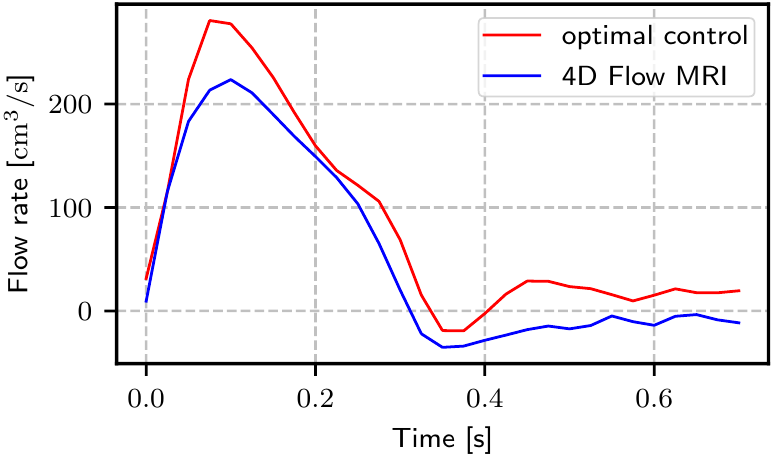}%
			\label{fig:4flow-dao}}
	}
	\caption{Comparison of 4D-Flow MRI flow waveforms and simulated flow waveforms with boundary conditions estimated with the proposed optimal control approach for case 4.}
	\label{fig:4flow-ocp}
\end{figure}

Additionally, Figure~\ref{fig:4flow-ocp} reports a comparison of time-dependent flow rates extracted from 4D-Flow MRI and those obtained from time-dependent simulations with boundary conditions estimated with optimal control. Even if the proposed estimation method assimilates only the average flow rates, the time-dependent waveforms are still recovered with a good degree of accuracy. For the sake of space, only results for case 4 are reported, with comparable results for the other cases. 

\begin{figure}[h]
	\centering
	\includegraphics[width=.9\linewidth]{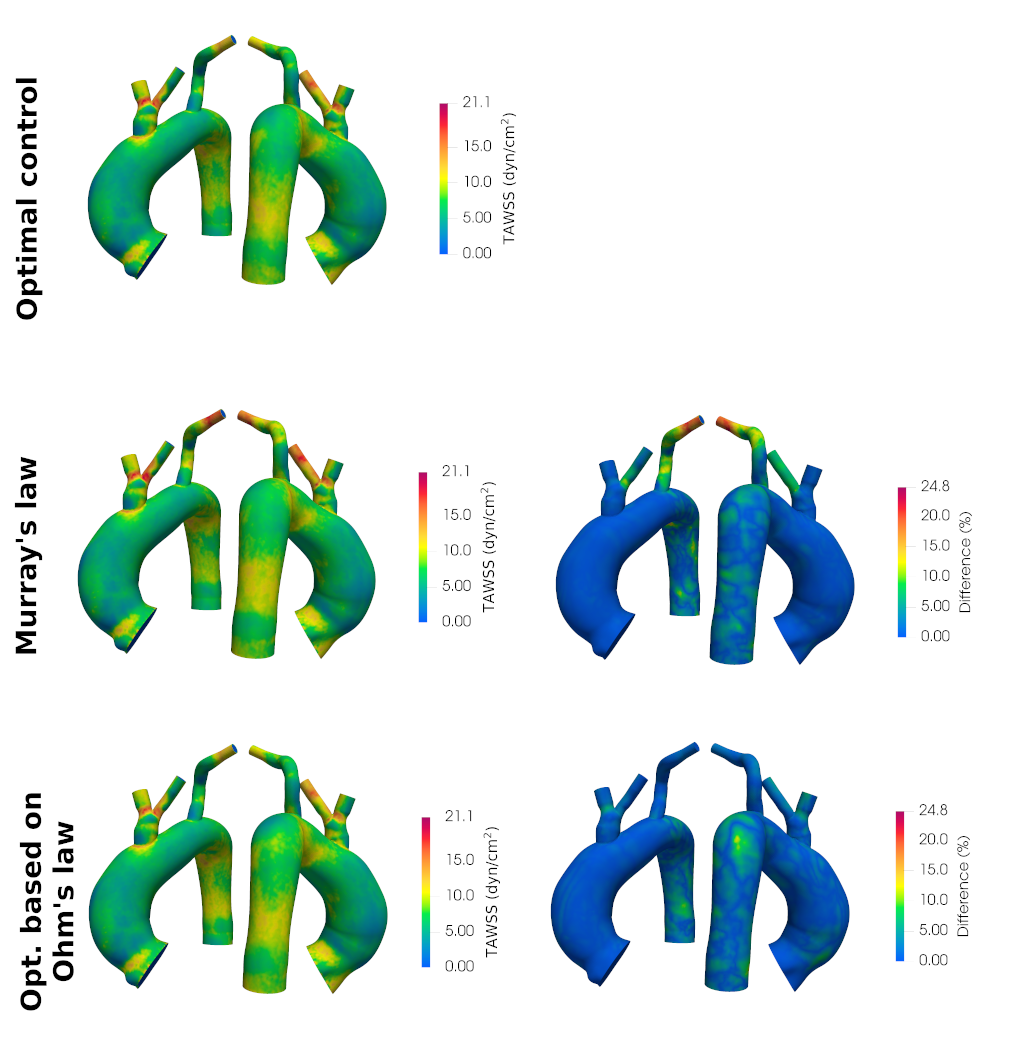}
	\caption{Left: Time-averaged wall shear stress (TAWSS) results obtained with the three different BC estimation techniques for case 3. Right: local absolute differences in TAWSS as compared to results obtained with optimal control. Anterior and posterior views of the anatomy are provided.}
	\label{fig:case010_tawss}
\end{figure}

\begin{figure}[h]
	\centering
	\includegraphics[width=.9\linewidth]{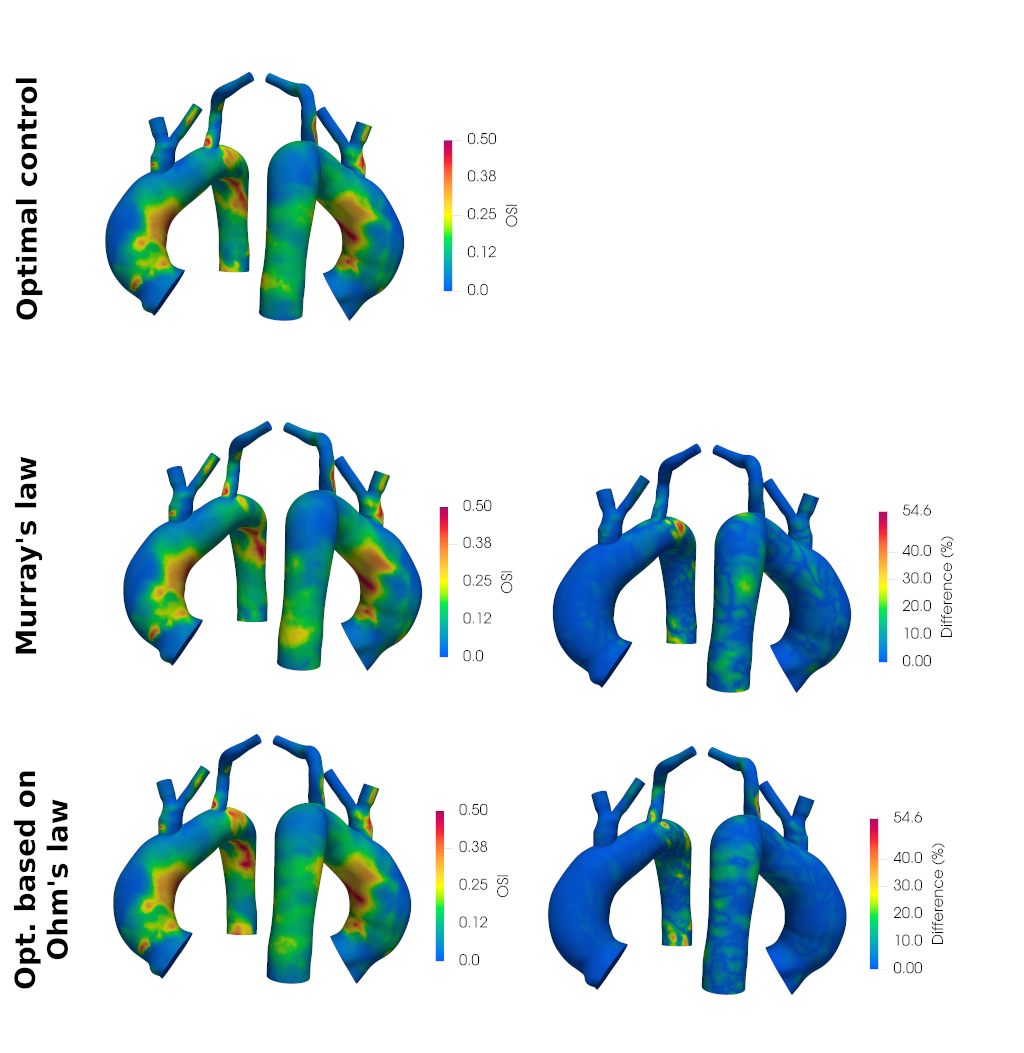}
	\caption{Left: Oscillatory Shear Index (OSI) results obtained with the three different BC estimation techniques for case 3. Right: local absolute differences in OSI as compared to results obtained with optimal control.
		Anterior and posterior views of the anatomy are provided.}
	\label{fig:case010_osi}
\end{figure}

To assess the influence that the adopted BCs estimation technique may have on clinically relevant parameters, we carried out an additional analysis on two relevant haemodynamic indicators, namely, the time-averaged wall shear stress (TAWSS) and the oscillatory shear index (OSI), calculated from Navier-Stokes simulation results using the equations reported by Martin et al.~\cite{martin2014computational}. 
For the sake of space, we report the analysis for case 3, but similar results were obtained for the other cases.
Fig.~\ref{fig:case010_tawss} shows, in the left column, the TAWSS obtained for the two cases with three different techniques (optimal control, Murray's law, and optimization based on Ohm's law), and in the right column the local relative difference with respect to optimal control results. The same analysis was repeated for the OSI in Fig.~\ref{fig:case010_osi}. For each point of the surface anatomy, the local relative difference was computed as:
\begin{equation}
\varepsilon_{r, TAWSS} = \frac{|\text{TAWSS} - \text{TAWSS}_{ocp}|}{max(\text{TAWSS}_{ocp})}
\end{equation}

\begin{equation}
\varepsilon_{r, OSI} = \frac{|\text{OSI} - \text{OSI}_{ocp}|}{max(\text{OSI}_{ocp})}
\end{equation}
For case 3, the difference in the TAWSS reaches a maximum relative difference of $24.8\%$, while the discrepancy in the OSI value reaches $55\%$.
The largest differences occur in the Murray's law case, in correspondence of LSUB, which is also where the estimated resistance values differ the most ($46\%$). 
It is worth noticing that, outside some specific `hot-spots', the relative error is generally lower, around 10\%. This analysis reveals the impact that the values of resistive boundary conditions have on haemodynamic indicators. In particular, given the relevance of TAWSS and OSI in a clinical context, the adoption of different techniques for BCs estimation could possibly affect the observations done by medical doctors, reaffirming the importance of adopting an automated, reliable, and operator-independent technique for boundary condition estimation.

\clearpage
\section{Discussion and Limitations} \label{section:limitations}
The proposed approach for boundary condition estimation presents some considerable advantages. 
First, optimal control relies on a rigorous and general mathematical framework, which ensures the possibility to apply it to different vessels in the cardiovascular system, provided that in-vivo measurements are available. 
Second, the mathematical structure of optimal control leaves great flexibility in terms of the nature and quantity of measurements to assimilate, an appealing feature for cardiovascular applications, where in-vivo measurements can be either pressure or flow waveforms. In particular, the proposed method does not specifically require data from 4D-Flow MRI, but it is compatible with any method providing flow rate information at the inlets and outlets of the region of interest, including the more common PC-MRI modality. 
Moreover, the proposed formulation assimilates data in a \textit{least-squares}  sense~\cite{gunzburger2002perspectives}, reducing the influence of stochastic measurement uncertainty on the solution.

On the other hand, the proposed technique presents some limitations.
The main limitation comes from the use of steady-state, linearized Navier-Stokes equations inside the optimal control problem. As explained in the Introduction, optimal control problems are characterized by a high computational cost. For this reason, the solution of a time-dependent optimal control problem is at the moment computationally intractable. The necessity to use a steady-state formulation of the problem requires resorting to the Stokes model, as Navier-Stokes equations may not converge to a steady-state solution for a complex flow like the one in the aortic arch. On one side, using a simplified model like the Stokes one significantly reduces the computational cost for the data assimilation process. On the other side, steady Stokes equations are clearly inadequate to get a realistic representation of the flow in the aorta. For this reason, the solution of the optimal control problem is used solely to estimate the outlet resistances, while for realistic pressure and velocity distributions a subsequent time-dependent simulation is still necessary.
One main consequence of using a steady-state formulation is the possibility to estimate only the total value of the resistance. Estimating the outlet capacitance, in fact, would require solving a time-dependent Navier-Stokes problem. For the same reason, it is not possible to determine multiple resistance values at each outlet when using more advanced outlet boundary conditions, such as the proximal and distal resistances in a 3-element Windkessel model, being then forced to adopt capacitance values and rules for resistance splitting taken from literature.
Lastly, the proposed method can only assimilate average flow rates, and not time-dependent flow waveforms. The results presented in Section~\ref{section:results}, in particular the simulated waveforms reported in Figure~\ref{fig:4flow-ocp}, prove that, when the resistance and capacitance values determined with the proposed method are used in an unsteady simulation, the obtained pressure and flow rates values are in sufficient agreement with the in-vivo measurements used for BC estimation.

\section{Summary and Conclusion} \label{section:conclusion}
In this work, we proposed a framework based on optimal control for the automated estimation of unknown resistance-type boundary conditions, while assimilating in-vivo pressure and flow rate measurements. 
The experiments conducted on four patient anatomies revealed the validity of the proposed optimal control-based technique for assimilating 4D-MRI data. 
Specifically, when compared to two other common techniques, namely, Murray's law and Ohm's law, the proposed framework performed consistently better.
An additional investigation of the effects of the different estimation methods on wall shear stress-related parameters exposed the influence that boundary conditions play on clinically relevant quantities. This further confirms the need for an automated method, such as the one proposed, which eliminates the expensive manual tuning phase, together with intra- and inter-operator variability. 
The proposed method, therefore, represents a first step in incorporating optimal control, thus a reliable, automated, and robust optimization technique into a framework which can be exploited by the medical community in a clinical setting. The field is promising and opens important perspectives for mathematical modelling and numerical simulation in cardiovascular flows.

A first extension of the present work is the application to other clinically relevant scenarios, such as coronary artery bypass grafts, where the scarcity and noise of available data limit the applicability of other estimation techniques.
Furthermore, the current framework estimates only resistive-type quantities, while the capacitance values are still chosen based on generic information available in the literature, usually combined with a manual tuning process. Future works will be directed toward the automated tuning of capacitance values as well.

\section*{Acknowledgements}
We acknowledge the support provided by the Canada Research Chairs program, the Radiological Society of North America, Compute Canada, and the H2020 ERC CoG AROMA-CFD project (G.A. 681447).

\bibliographystyle{plain}
\bibliography{bibliography}

\end{document}